\documentclass[twoside,11pt]{article}

\usepackage{jmlr2e}

\usepackage{amsfonts}
\usepackage{epsfig}
\usepackage{amsmath}
\usepackage{amssymb}
\usepackage{algorithm}
\usepackage{algorithmic}
\usepackage{subfig}
\usepackage{graphicx}
\usepackage{color}
\usepackage[section] {placeins}
\usepackage{epstopdf}
\usepackage{url}

\usepackage{natbib}
\usepackage{hyperref}

\newcounter{ass_counter}
\newtheorem{assumption}[ass_counter]{Assumption}

\def\E{\mathbb{E}}
\def\PO{\mathcal{P}_{\Omega}}
\def\PS{\mathcal{P}_{S}}
\def\PO{\mathcal{P}_{\Omega}}

\newcommand{\SynGD}{\mbox{\sc SynGD}}
\newcommand{\ASCD}{\mbox{\sc AsySCD}}
\newcommand{\AsySCD}{\mbox{\sc AsySCD}}
\newcommand{\Lres}{L_{\mbox{\rm\scriptsize res}}}
\newcommand{\Lmax}{L_{\max}}
\newcommand{\beq}{\begin{equation}}
\newcommand{\eeq}{\end{equation}}
\def\eqnok#1{(\ref{#1})}
\newcommand{\R}{\mathbb{R}}

\def\hogwild{{\sc Hogwild!}}

\jmlrheading{x}{xxxx}{xx-xx}{xx/xx}{xx/xx}{Liu, Wright, R{\'e}, Bittorf, Sridhar}

\ShortHeadings{$\AsySCD$}{Liu, Wright, R{\'e}, Bittorf, Sridhar}
\firstpageno{1}

\begin{document}

\title{An Asynchronous Parallel Stochastic Coordinate Descent Algorithm}

\author{\name Ji Liu \email ji.liu.uwisc@gmail.com \\
       \name Stephen J. Wright \email swright@cs.wisc.edu \\
              \addr Department of Computer Sciences, University of Wisconsin-Madison, Madison, WI 53706 \\
       \name Christopher R{\'e} \email chrismre@cs.stanford.edu \\
       \addr Department of Computer Science, Stanford University, Stanford, CA 94305 \\
       \name Victor Bittorf \email victor.bittorf@cloudera.com \\
       \addr Cloudera, Inc., Palo Alto, CA 94304 \\
       \name Srikrishna Sridhar \email srikris@graphlab.com \\
       \addr GraphLab Inc., Seattle, WA  98103 \\
       }
      
     
\editor{Leon Bottou}

\maketitle


%




\begin{abstract}
We describe an asynchronous parallel stochastic coordinate descent
algorithm for minimizing smooth unconstrained or separably constrained
functions. The method achieves a linear convergence rate on functions
that satisfy an essential strong convexity property and a sublinear
rate ($1/K$) on general convex functions. Near-linear speedup on a
multicore system can be expected if the number of processors is
$O(n^{1/2})$ in unconstrained optimization and $O(n^{1/4})$ in the
separable-constrained case, where $n$ is the number of variables. We
describe results from implementation on 40-core processors.
\end{abstract}

\begin{keywords}
Asynchronous Parallel Optimization, Stochastic Coordinate Descent
\end{keywords}


\section{Introduction} \label{sec:intro}

Consider the convex optimization problem
\begin{equation}
\label{eqn_mainproblem}
\min_{x \in \Omega} \, \quad f(x),
\end{equation}
where $\Omega \subset \mathbb{R}^{n}$ is a closed convex set and $f$
is a smooth convex mapping from an open neighborhood of $\Omega$ to
$\mathbb{R}$. We consider two particular cases of $\Omega$ in this
paper: the unconstrained case $\Omega= \mathbb{R}^n$, and the
separable case
\beq \label{eq:sep}
\Omega = \Omega_1 \times \Omega_2 \times \dotsc \times
\Omega_n,
\eeq
where each $\Omega_i$, $i=1,2,\dotsc,n$ is a closed subinterval of the
real line.

Formulations of the type (\ref{eqn_mainproblem},\ref{eq:sep}) arise in
many data analysis and machine learning problems, for example, support
vector machines (linear or nonlinear dual formulation)
\citep{CortesVapnik95}, LASSO (after decomposing $x$ into positive and
negative parts) \citep{Tibshirani96LASSO}, and logistic
regression. Algorithms based on gradient and approximate or partial
gradient information have proved effective in these settings. We
mention in particular gradient projection and its accelerated variants
\citep{nesterov2004introductory}, accelerated proximal gradient
methods for regularized objectives \citep{BeckT09}, and stochastic
gradient methods \citep{Nemirovski09, Shamir2013icml}. These methods
are inherently serial, in that each iteration depends on the result of
the previous iteration. Recently, parallel multicore versions of
stochastic gradient and stochastic coordinate descent have been
described for problems involving large data sets; see for example
\citet{Hogwild11nips, Richtarik12arXiv, Avron13arXiv}.

This paper proposes an asynchronous stochastic coordinate descent
($\ASCD$) algorithm for convex optimization. Each step of $\ASCD$
chooses an index $i \in \{1,2,\dotsc,n\}$ and subtracts a short,
constant, positive multiple of the $i$th partial gradient $\nabla_i
f(x) := \partial f / \partial x_i$ from the $i$th component of
$x$. When separable constraints \eqnok{eq:sep} are present, the update
is ``clipped'' to maintain feasibility with respect to
$\Omega_i$. Updates take place in parallel across the cores of a
multicore system, without any attempt to synchronize computation
between cores.
We assume that there is a bound $\tau$ on the age of the updates, that
is, no more than $\tau$ updates to $x$ occur between the time at which
a processor reads $x$ (and uses it to evaluate one element of the
gradient) and the time at which this processor makes its update to a
single element of $x$.  (A similar model of parallel asynchronous
computation was used in \hogwild~\citep{Hogwild11nips}.)  Our
implementation, described in Section~\ref{sec_exp}, is a little more
complex than this simple model would suggest, as it is tailored to the
architecture of the Intel Xeon machine that we use for experiments.

We show that linear convergence can be attained if an ``essential
strong convexity'' property \eqnok{eq:esc} holds, while sublinear
convergence at a ``$1/K$'' rate can be proved for general convex
functions. Our analysis also defines a sufficient condition for
near-linear speedup in the number of cores used. This condition
relates the value of delay parameter $\tau$ (which relates to the
number of cores / threads used in the computation) to the problem
dimension $n$. A parameter that quantifies the cross-coordinate
interactions in $\nabla f$ also appears in this relationship. When the
Hessian of $f$ is nearly diagonal, the minimization problem can almost
be separated along the coordinate axes, so higher degrees of
parallelism are possible.


We review related work in
Section~\ref{sec_relatedwork}. Section~\ref{sec_alg} specifies the
proposed algorithm.  Convergence results for unconstrained and
constrained cases are described in Sections~\ref{sec_unconstrained}
and \ref{sec_constrained}, respectively, with proofs given in the
appendix. Computational experience is reported in
Section~\ref{sec_exp}. We discuss several variants of $\ASCD$ in
Section~\ref{sec_extension}.  Some conclusions are given in
Section~\ref{sec_conclusion}.

\subsection*{Notation and Assumption}\label{sec_NA}
We use the following notation.
\begin{itemize}
\item $e_i\in\R^{n}$ denotes the $i$th natural basis vector $(0,
  \dotsc,0,1,0,\dotsc,0)^T$ with the `''$1$'' in the $i$th position.
\item $\| \cdot \|$ denotes the Euclidean norm $\|\cdot \|_2$.
\item $S \subset \Omega$ denotes the set on which $f$ attains its
  optimal value, which is denoted by $f^*$.
\item $\PS(\cdot)$ and $\PO(\cdot)$ denote Euclidean projection onto
  $S$ and $\Omega$, respectively.
\item We use $x_i$ for the $i$th element of $x$, and $\nabla_i f(x)$
for the $i$th element of the gradient vector $\nabla f(x)$.
\item We define the following {\em essential strong convexity}
  condition for a convex function $f$ with respect to the optimal set
  $S$, with parameter $l>0$:
\begin{equation}
f(x) - f(y) \geq \langle \nabla f(y), x-y \rangle + {l\over 2}\|x-y\|^2 \quad
\mbox{for all $x, y \in \Omega$ with $\PS(x)=\PS(y)$.} \label{eq:esc}
\end{equation}
%
This condition is significantly weaker than the usual strong convexity
condition, which requires the inequality to hold for {\em all} $x, y
\in \Omega$. In particular, it allows for non-singleton solution sets
$S$, provided that $f$ increases at a uniformly quadratic rate with
distance from $S$.  (This property is noted for convex quadratic $f$
in which the Hessian is rank deficient.)  Other examples of
essentially strongly convex functions that are not strongly convex
include:
\begin{itemize}
\item $f(Ax)$ with arbitrary linear transformation $A$, where
  $f(\cdot)$ is strongly convex;
\item $f(x) = \max(a^Tx - b, 0)^2$, for $a \neq 0$.
\end{itemize}

\item Define $\Lres$ as the {\em restricted Lipschitz constant} for
  $\nabla f$, where the ``restriction'' is to the coordinate
  directions: We have
\[
\| \nabla f(x) - \nabla f(x+te_i)\| \leq \Lres |t|, \quad
\mbox{for all $i=1,2,\dotsc,n$ and $t \in \R$, with $x, x+t e_i \in \Omega$.}
\]
\item Define $L_i$ as the {\em coordinate Lipschitz constant} for
  $\nabla f$ in the $i$th coordinate direction: We have
\[
f(x+te_i) - f(x) \leq \langle \nabla_if(x),~t \rangle + {L_i\over 2}t^2,
\quad \mbox{for $i \in \{1,2,\dotsc,n\}$, and $x,x+t e_i \in \Omega$,}
\]
or equivalently
\[
|\nabla_if(x) - \nabla_if(x+te_i)| \leq L_i|t|.
\]
\item $\Lmax:=\max_{i=1,2,\dotsc,n} L_i$. 
\item For the initial point $x_0$, we define
\begin{equation} \label{eq:R0.def}
R_0 := \| x_0 - \PS(x_0) \|.
\end{equation}
\end{itemize}
Note that $\Lres\geq \Lmax$.

We use $\{ x_j \}_{j=0,1,2,\dotsc}$ to denote the sequence of iterates
generated by the algorithm from starting point $x_0$. Throughout the
paper, we assume that $S$ is nonempty.


\subsection*{Lipschitz Constants}

The nonstandard Lipschitz constants $\Lres$, $\Lmax$, and $L_i$,
$i=1,2,\dotsc,n$ defined above are crucial in the analysis of our
method. Besides bounding the nonlinearity of $f$ along various
directions, these quantities capture the interactions between the
various components in the gradient $\nabla f$, as quantified in the
off-diagonal terms of the Hessian $\nabla^2 f(x)$ --- although the
stated conditions do not require this matrix to exist.

We have noted already that $\Lres/\Lmax \ge 1$. Let us consider upper
bounds on this ratio under certain conditions. When $f$ is twice
continuously differentiable, we have
\[
L_{\max} = \sup_{x \in \Omega} \, \max_{i=1,2,\dotsc,n} \, [\nabla^2
  f(x)]_{ii}.
\]
Since $\nabla^2 f (x) \succeq 0$ for $x \in \Omega$, we have that
\[
| [\nabla^2 f(x)]_{ij} | \le \sqrt{L_i L_j} \le \Lmax, \quad \forall \, i,j=1,2,\dotsc,n.
\]
Thus $\Lres$, which is a bound on the largest column norm for
$\nabla^2 f(x)$ over all $x \in \Omega$, is bounded by $\sqrt{n}
\Lmax$, so that
\[
\frac{\Lres}{\Lmax} \le \sqrt{n}.
\]
If the Hessian is structurally sparse, having at most $p$ nonzeros per
row/column, the same argument leads to $\Lres/\Lmax \le \sqrt{p}$.

If $f(x)$ is a convex quadratic with Hessian $Q$, we have
\[
\Lmax = \max_i \, Q_{ii}, \quad
\Lres = \max_i \| Q_{\cdot i} \|_2,
\]
where $Q_{\cdot i}$ denotes the $i$th column of $Q$. If $Q$ is
diagonally dominant, we have for any column $i$ that
\[
\| Q_{\cdot i} \|_2 \le Q_{ii} + \| [Q_{ji}]_{j \neq i} \|_2 \le
Q_{ii} + \sum_{j \neq i} |Q_{ji}| \le 2 Q_{ii},
\]
which, by taking the maximum of both sides, implies that $\Lres/\Lmax
\le 2$ in this case.

Finally, consider the objective $f(x)=\frac12 \|Ax-b\|^2$ and assume
that $A\in\mathbb{R}^{m\times n}$ is a random matrix whose entries are
i.i.d from $\mathcal{N}(0, 1)$.  The diagonals of the Hessian are
$A_{\cdot i}^T A_{\cdot i}$ (where $A_{\cdot i}$ is the $i$th column
of $A$), which have expected value $m$, so we can expect $\Lmax$ to be
not less than $m$.  Recalling that $\Lres$ is the maximum column norm of
$A^TA$, we have
\begin{align*}
\E (\|A^TA_{\cdot i}\|) & \leq \E (|A_{\cdot i}^TA_{\cdot i} |) +
\E (\| [ A_{\cdot j}^T A_{\cdot i}]_{j \neq i} \|) \\
& = m +  \E \sqrt{\sum_{j\neq i}|A_{\cdot j}^TA_{\cdot i} |^2} \\
& \leq m + \sqrt{\sum_{j\neq i}\E|A_{\cdot j}^T A_{\cdot i} |^2} \\
& = m + \sqrt{(n-1)m},
\end{align*}
where the second inequality uses Jensen's inequality
and the final equality uses
\[
\E (|A_{\cdot j}^TA_{\cdot i} |^2) = \E(A_{\cdot j}^T\E(A_{\cdot
  i}A_{\cdot i}^T) A_{\cdot j})
  =\E(A_{\cdot j}^TI A_{\cdot
  j})=\E(A_{\cdot j}^TA_{\cdot j})=m.
\]
We can thus estimate the upper bound on $\Lres/\Lmax$ roughly by $1+
\sqrt{n/m}$ for this case. 


\section{Related Work}\label{sec_relatedwork}



This section reviews some related work on coordinate relaxation and
stochastic gradient algorithms.


Among {\em cyclic coordinate descent} algorithms, \citet{Tseng01}
proved the convergence of a block coordinate descent method for
nondifferentiable functions with certain conditions. Local and global
linear convergence were established under additional assumptions, by
\citet{LuoTseng92} and \citet{WangLin13}, respectively. Global linear
(sublinear) convergence rate for strongly (weakly) convex optimization
was proved by \citet{Beck13}. Block-coordinate approaches based on
proximal-linear subproblems are described by
\citet{TseY06,TseY07a}. \citet{Wright12} uses acceleration on reduced
spaces (corresponding to the optimal manifold) to improve the local
convergence properties of this approach.


{\em Stochastic coordinate descent} is almost identical to cyclic
coordinate descent except selecting coordinates in a random
manner. \citet{Nesterov12} studied the convergence rate for a
stochastic block coordinate descent method for unconstrained and
separably constrained convex smooth optimization, proving linear
convergence for the strongly convex case and a sublinear $1/K$ rate
for the convex case. Extensions to minimization of composite functions
are described by \citet{Richtarik11} and \citet{LuXiao13}.

{\em Synchronous parallel methods} distribute the workload and data
among multiple processors, and coordinate the computation among
processors. \citet{Ferris94} proposed to distribute variables among
multiple processors and optimize concurrently over each subset. The
synchronization step searches the affine hull formed by the current
iterate and the points found by each processor. Similar ideas appeared
in \citep{Mangasarian95}, with a different synchronization
step. \citet{Ma12} considered a multiple splitting algorithm for
functions of the form $f(x) = \sum_{k=1}^N f_k(x)$ in which $N$ models
are optimized separately and concurrently, then combined in an
synchronization step. The alternating direction method-of-multiplier
(ADMM) framework \citep{Boyd11} can also be implemented in parallel.
This approach dissects the problem into multiple subproblems (possibly after
replication of primal variables) and optimizes concurrently, then
synchronizes to update multiplier estimates. \citet{Duchi12} described
a subgradient dual-averaging algorithm for partially separable
objectives, with subgradient evaluations distributed between cores and
combined in ways that reflect the structure of the objective. Parallel
stochastic gradient approaches have received broad attention; see
\citet{AgarwalD12} for an approach that allows delays between
evaluation and update, and \citet{Cotter11} for a minibatch stochastic
gradient approach with Nesterov
acceleration. \citet{Shalev-Shwartz2013} proposed an accelerated
stochastic dual coordinate ascent method.


 Among {\em synchronous parallel methods for (block) coordinate
   descent}, \citet{Richtarik12arXiv} described a method of this type
 for convex composite optimization problems.  All processors update
 randomly selected coordinates or blocks, concurrently and
 synchronously, at each iteration. Speedup depends on the sparsity of
 the data matrix that defines the loss functions. Several variants
 that select blocks greedily are considered by \citet{ScherrerTHH12}
 and \citet{Yin13}. \citet{Yang13} studied the parallel stochastic dual 
 coordinate ascent method and emphasized the balance between 
 computation and communication. 


We turn now to {\em asynchronous parallel
  methods}. \citet{Bertsekas89} introduced an asynchronous parallel
implementation for general fixed point problems $x = q(x)$ over a
separable convex closed feasible region. (The optimization problem
\eqnok{eqn_mainproblem} can be formulated in this way by defining
$q(x) := \mathcal{P}_\Omega[(I - \alpha \nabla f)(x)]$ for some fixed
$\alpha>0$.) Their analysis allows inconsistent reads for $x$, that
is, the coordinates of the read $x$ have different ``ages.''  Linear
convergence is established if all ages are bounded and $\nabla^2 f(x)$
satisfies a diagonal dominance condition guaranteeing that the
iteration $x=q(x)$ is a maximum-norm contraction mapping for
sufficient small $\alpha$. However, this condition is strong ---
stronger, in fact, than the strong convexity condition. For convex
quadratic optimization $f(x) = {1\over 2}x^TAx + bx$, the contraction
condition requires diagonal dominance of the Hessian: $A_{ii} >
\sum_{i\neq j}|A_{ij}|$ for all $i=1,2,\dotsc,n$. By comparison,
$\ASCD$ guarantees linear convergence rate under the essential strong
convexity condition \eqnok{eq:esc}, though we do not allow
inconsistent read. (We require the vector $x$ used for each evaluation
of $\nabla_i f(x)$ to have existed at a certain point in time.)


\hogwild~\citep{Hogwild11nips} is a lock-free, asynchronous parallel
implementation of a stochastic-gradient method, targeted to a
multicore computational model similar to the one considered here. Its
analysis assumes consistent reading of $x$, and it is implemented
without locking or coordination between processors. Under certain
conditions, convergence of \hogwild~approximately matches the
sublinear $1/K$ rate of its serial counterpart, which is the
constant-steplength stochastic gradient method analyzed in
\citet{Nemirovski09}.


We also note recent work by \citet{Avron13arXiv}, who proposed an
asynchronous linear solver to solve $Ax=b$ where $A$ is a symmetric
positive definite matrix, proving a linear convergence rate. Both
inconsistent- and consistent-read cases are analyzed in this paper,
with the convergence result for inconsistent read being slightly
weaker.

The $\AsySCD$ algorithm described in this paper was extended to solve
the composite objective function consisting of a smooth convex
function plus a separable convex function in a later work
\citep{liu2014asynchronous}, which pays particular attention to the
inconsistent-read case.


\section{Algorithm}\label{sec_alg}

In $\ASCD$, multiple processors have access to a shared data structure
for the vector $x$, and each processor is able to compute a randomly
chosen element of the gradient vector $\nabla f(x)$. Each processor
repeatedly runs the following coordinate descent process (the
steplength parameter $\gamma$ is discussed further in the next
section):
\begin{itemize}
\item[R:] Choose an index $i \in \{1,2,\dotsc,n\}$ at random, read $x$,
  and evaluate $\nabla_i f(x)$;
\item[U:] Update component $i$ of the shared $x$ by taking a step of
  length $\gamma/\Lmax$ in the direction $-\nabla_i f(x)$.
\end{itemize}

\begin{algorithm}[ht!]
\caption{Asynchronous Stochastic Coordinate Descent Algorithm $x_{K+1}=\ASCD(x_0, \gamma, K)$}
\label{alg_ascd}
\begin{algorithmic}[1]
\REQUIRE $x_0\in\Omega$, $\gamma$, and $K$
\ENSURE $x_{K+1}$
\STATE Initialize $j \leftarrow 0$;
\WHILE{$j \leq K$}
\STATE Choose $i(j)$ from $\{1,\dotsc,n\}$ with equal probability;
\STATE $x_{j+1} \leftarrow \mathcal{P}_{\Omega}\left(x_j-\frac{\gamma}{\Lmax} e_{i(j)} \nabla_{i(j)}f(x_{k(j)})\right)$; \label{step_proj}
\STATE $j \leftarrow  j +1$;
\ENDWHILE
\end{algorithmic}
\end{algorithm}

Since these processors are being run concurrently and without
synchronization, $x$ may change between the time at which it is read
(in step R) and the time at which it is updated (step U). We capture
the system-wide behavior of $\ASCD$ in Algorithm~\ref{alg_ascd}.
There is a global counter $j$ for the total number of updates; $x_j$
denotes the state of $x$ after $j$ updates. The index $i(j) \in
\{1,2,\dotsc,n\}$ denotes the component updated at step $j$. $k(j)$
denotes the $x$-iterate at which the update applied at iteration $j$
was calculated. Obviously, we have $k(j) \le j$, but we assume that
the delay between the time of evaluation and updating is bounded
uniformly by a positive integer $\tau$, that is, $j-k(j) \le \tau$ for
all $j$. The value of $\tau$ captures the essential parallelism in the
method, as it indicates the number of processors that are involved in
the computation.

The projection operation $P_{\Omega}$ onto the feasible set is not
needed in the case of unconstrained optimization. For separable
constraints \eqnok{eq:sep}, it requires a simple clipping operation on
the $i(j)$ component of $x$.

We note several differences with earlier asynchronous approaches.
Unlike the asynchronous scheme in \citet[Section 6.1]{Bertsekas89},
the {\em latest} value of $x$ is updated at each step, not an earlier
iterate.  Although our model of computation is similar to
\hogwild~\citep{Hogwild11nips}, the algorithm differs in that each
iteration of $\ASCD$ evaluates a single component of the gradient
exactly, while \hogwild~ computes only a (usually crude) estimate of
the full gradient. Our analysis of $\ASCD$ below is comprehensively
different from that of \citet{Hogwild11nips}, and we obtain stronger
convergence results.

\section{Unconstrained Smooth Convex Case} \label{sec_unconstrained}

This section presents results about convergence of $\ASCD$ in the
unconstrained case $\Omega = \R^n$. The theorem encompasses both the
linear rate for essentially strongly convex $f$ and the sublinear rate
for general convex $f$. The result depends strongly on the delay
parameter $\tau$. (Proofs of results in this section appear in
Appendix~\ref{app:unc}.) In Algorithm~\ref{alg_ascd}, the indices
$i(j)$, $j=0,1,2,\dotsc$ are random variables. We denote the
expectation over all random variables as $\E$, the conditional
expectation in term of $i(j)$ given $i(0), i(1), \cdots, i(j-1)$ as
$\E_{i(j)}$.

A crucial issue in $\ASCD$ is the choice of steplength parameter
$\gamma$. This choice involves a tradeoff: We would like $\gamma$ to
be long enough that significant progress is made at each step, but not
so long that the gradient information computed at step $k(j)$ is stale
and irrelevant by the time the update is applied at step $j$. We
enforce this tradeoff by means of a bound on the ratio of expected
squared norms on $\nabla f$ at successive iterates; specifically,
\begin{equation} \label{eq:ratio_unc}
\rho^{-1}\leq\frac{\E\|\nabla f(x_{j+1})\|^2}{\E\|\nabla f(x_{j})\|^2}\leq \rho,
\end{equation}
where $\rho > 1$ is a user defined parameter. 
The analysis becomes a delicate balancing act in the choice of $\rho$
and steplength $\gamma$ between aggression and excessive
conservatism. We find, however, that these values can be chosen to
ensure steady convergence for the asynchronous method at a {\em
  linear} rate, with rate constants that are almost consistent with
vanilla short-step full-gradient descent.

We use the following assumption in some of the results of this
section.
\begin{assumption} \label{ass_1}
There is a real number $R$ such that
\[
\| x_j - \PS(x_j) \| \le R, \quad \mbox{for all $j=0,1,2,\dotsc$.}
\]
\end{assumption}
Note that this assumption is not needed in our convergence results in
the case of strongly convex functions.  in our theorems below, it is
invoked only when considering {\em general} convex functions.


%
\begin{theorem}  \label{AsySCD:thm_1}
Suppose that $\Omega=\R^n$ in \eqnok{eqn_mainproblem}. For any
$\rho>1$, define the quantity $\psi$ as follows:
\beq \label{eq:defpsi}
\psi := 1+ \frac{2 \tau \rho^{\tau} \Lres}{\sqrt{n} \Lmax}.
\eeq
Suppose that the steplength parameter $\gamma>0$ satisfies the following
three upper bounds:
\begin{subequations}
\label{eq:boundgamma}
\begin{align}
\label{eq:boundgamma.1}
\gamma & \le \frac{1}{\psi}, \\
\label{eq:boundgamma.2}
\gamma & \le \frac{(\rho-1) \sqrt{n} \Lmax}{2 \rho^{\tau+1} \Lres}, \\
\label{eq:boundgamma.3}
\gamma & \le \frac{(\rho-1) \sqrt{n} \Lmax}{\Lres \rho^{\tau} (2+\frac{\Lres}{\sqrt{n} \Lmax})}.
\end{align}
\end{subequations}
Then we have that for any $j\geq 0$ that
\begin{equation}
\rho^{-1}\E(\|\nabla f(x_j)\|^2)  \leq {\E (\|\nabla f(x_{j+1})\|^2)}\leq \rho \E(\|\nabla f(x_j)\|^2).
\label{eqn_thm_1}
\end{equation}
Moreover, if the essentially strong convexity property \eqnok{eq:esc}
holds with $l>0$, we have
\beq
\E (f(x_{j}) - f^*) \leq \left(1 - \frac{2l\gamma}{n\Lmax}
\left( 1-\frac{\psi}{2} \gamma \right)\right)^{j} (f(x_{0})-f^*). \label{eqn_thm_2}
\eeq
For general smooth convex functions $f$, assuming additionally
that Assumption~\ref{ass_1} holds, we have
\beq
\E (f(x_j) - f^*) \leq \frac{1}{(f(x_{0}) -f^*)^{-1} +
j \gamma (1-\frac{\psi}{2} \gamma) /(n \Lmax R^2) }. \label{eqn_thm_3}
\eeq
\end{theorem}
This theorem demonstrates linear convergence \eqnok{eqn_thm_2} for
$\ASCD$ in the unconstrained essentially strongly convex case. This
result is better than that obtained for
\hogwild~\citep{Hogwild11nips}, which guarantees only sublinear
convergence under the stronger assumption of strict convexity.

The following corollary proposes an interesting particular choice of
the parameters for which the convergence expressions become more
comprehensible. The result requires a condition on the delay bound
$\tau$ in terms of $n$ and the ratio $\Lmax/\Lres$.
\begin{corollary} \label{co:thm_1}
Suppose that 
\beq \label{eq:boundtau}
\tau + 1 \leq \frac{\sqrt{n} \Lmax}{2e\Lres}.
\eeq
Then if we choose
\beq \label{eq:choicerho}
\rho=1+{2e\Lres\over  \sqrt{n}\Lmax},
\eeq
define $\psi$ by \eqnok{eq:defpsi}, and set $\gamma=1/\psi$, 
we have for the essentially
strongly convex case \eqnok{eq:esc} with $l>0$ that
\beq \label{eqn_thm_2_good}
\E (f(x_j)-f^*) \le \left( 1-\frac{l}{2n\Lmax} \right)^j (f(x_0)-f^*).
\eeq
For the case of general convex $f$, if we assume additionally
that Assumption~\ref{ass_1} is satisfied, we have
\beq \label{eqn_thm_3_good}
\E (f(x_j)-f^*) \le  \frac{1}{(f(x_0)-f^*)^{-1} + {j}/(4n \Lmax R^2)}.
\eeq
\end{corollary}


We note that the linear rate \eqnok{eqn_thm_2_good} is broadly
consistent with the linear rate for the classical steepest descent
method applied to strongly convex functions, which has a rate constant
of $(1-2l/L)$, where $L$ is the standard Lipschitz constant for
$\nabla f$. If we assume (not unreasonably) that $n$ steps of
stochastic coordinate descent cost roughly the same as one step of
steepest descent, and note from \eqnok{eqn_thm_2_good} that $n$ steps
of stochastic coordinate descent would achieve a reduction factor of
about $(1-l/(2 \Lmax))$, a standard argument would suggest that
stochastic coordinate descent would require about $4 \Lmax/L$ times
more computation. (Note that $\Lmax/L \in [1/n,1]$.)
The stochastic approach may gain an advantage from the parallel
implementation, however. Steepest descent  requires
synchronization and careful division of gradient evaluations, whereas the
stochastic approach can be implemented in an asynchronous fashion.

For the general convex case, \eqnok{eqn_thm_3_good} defines a
sublinear rate, whose relationship with the rate of the steepest
descent for general convex optimization is similar to the previous
paragraph. 

As noted in Section~\ref{sec:intro}, the parameter $\tau$ is closely
related to the number of cores that can be involved in the
computation, without degrading the convergence performance of the
algorithm. In other words, if the number of cores is small enough such
that \eqref{eq:boundtau} holds, the convergence expressions
\eqref{eqn_thm_2_good}, \eqref{eqn_thm_3_good} do not depend on the
number of cores, implying that linear speedup can be expected. A small
value for the ratio $\Lres / \Lmax$ (not much greater than $1$)
implies a greater degree of potential parallelism. As we note at the
end of Section~\ref{sec:intro}, this ratio tends to be small in some
important applications --- a situation that would allow $O(\sqrt{n})$
cores to be used with near-linear speedup.


We conclude this section with a high-probability estimate for
convergence of  the sequence of function values.
%
\begin{theorem} \label{thm:co1}
Suppose that the assumptions of Corollary~\ref{co:thm_1} hold,
including the definitions of $\rho$ and $\psi$. Then for any
$\epsilon \in (0,f(x_0)-f^*)$
and $\eta\in(0,1)$, we have that
\beq \label{eq:thm:co1_2}
\mathbb{P}\left(f(x_j)-f^*\le \epsilon \right) \ge 1-\eta,
\eeq
provided that either of the following sufficient conditions hold for
the index $j$. In the essentially strongly convex case~\eqnok{eq:esc}
with $l>0$, it suffices to have
\beq \label{eq:def_j1}
j \geq \frac{2n\Lmax}{l} \left( {\log{f(x_0)-f^*\over \epsilon\eta}}\right).
\eeq
For the general convex case, if we assume additionally that
Assumption~\ref{ass_1} holds, a sufficient condition is
\beq \label{eq:def_j2}
j \geq 4n\Lmax R^2\left({1\over \epsilon\eta} - {1\over f(x_0)-f^*}\right).
\eeq
\end{theorem}


\section{Constrained Smooth Convex Case}\label{sec_constrained}


This section considers the case of separable constraints
\eqnok{eq:sep}. We show results about convergence rates and
high-probability complexity estimates, analogous to those of the
previous section.  Proofs appear in Appendix~\ref{app:con}.



As in the unconstrained case, the steplength $\gamma$ should be chosen
to ensure steady progress while ensuring that update information does
not become too stale. Because constraints are present, the ratio
\eqnok{eq:ratio_unc} is no longer appropriate. We use instead a ratio
of squares of expected differences in successive primal iterates:
\begin{equation}
\frac{\E\|x_{j-1}-\bar{x}_j\|^2 }{\E \|x_j-\bar{x}_{j+1}\|^2}, \label{eq:cn_ratio}
\end{equation}
where $\bar{x}_{j+1}$ is the hypothesized full update obtained by
applying the single-component update to {\em every} component of
$x_j$, that is,
\begin{equation} \label{eq:def.xbar}
\bar{x}_{j+1} := \arg\min_{x\in\Omega} \, \langle \nabla f(x_{k(j)}),
x-x_j \rangle + \frac{\Lmax}{2\gamma}\|x-x_j\|^2.
\end{equation}
In the unconstrained case $\Omega=\R^n$, the ratio \eqref{eq:cn_ratio}
reduces to
\[
\frac{\E\|\nabla f(x_{k(j-1)})\|^2}{\E\|\nabla f(x_{k(j)})\|^2},
\]
which is evidently related to \eqnok{eq:ratio_unc}, but not identical.

We have the following result concerning convergence of the expected
error to zero.
\begin{theorem} \label{thm_2}
Suppose that $\Omega$ has the form \eqnok{eq:sep} and that $n \ge
5$. Let $\rho$ be a constant with
$\rho>\left(1-2/\sqrt{n}\right)^{-1}$, and define the quantity $\psi$
as follows:
\beq \label{eq:defpsic} \psi:=
     {1+\frac{\Lres\tau\rho^{\tau}}{\sqrt{n}\Lmax}\left({2
         +{\Lmax\over \sqrt{n}\Lres} + {2\tau\over n}}\right)}.  
\eeq
Suppose that the steplength parameter $\gamma>0$ satisfies the
following two upper bounds:
  \beq \label{eq:boundgammac}
\gamma \le {1\over \psi},\quad \gamma
  \le \left(1-{1\over \rho}-{2\over
    \sqrt{n}}\right)\frac{\sqrt{n}\Lmax}{4\Lres\tau\rho^{\tau}}.
\eeq
Then we have
\begin{equation}
\E\|x_{j-1}-\bar{x}_j\|^2 \leq \rho \E \|x_j-\bar{x}_{j+1}\|^2,
\quad j=1,2,\dotsc. \label{eqn_thm2_1}
\end{equation}
If the essential strong convexity property~\eqref{eq:esc}
holds with $l>0$, we have for $j=1,2,\dotsc$ that
\begin{align}  
&\E \|x_{j} - \PS(x_{j})\|^2 + \frac{2\gamma}{\Lmax}(\E f(x_{j}) - f^*) \label{eqn_thm2_3}
\\
&\leq \left(1- \frac{l} {n(l+\gamma^{-1}\Lmax)} \right)^{j} \left( R_0^2 + \frac{2\gamma}{\Lmax}(f(x_0) - f^*)\right) \nonumber 
\end{align}
where $R_0$ is defined in \eqref{eq:R0.def}.
For general smooth convex function $f$, we have
\beq  \label{eqn_thm2_2}
\E f(x_{j})- f^* \leq \frac{n(R_0^2\Lmax+ 2\gamma(f(x_0)- f^*))}{2\gamma (n+j)}.
\eeq
\end{theorem}

Similarly to the unconstrained case, the following corollary proposes
an interesting particular choice for the parameters for which the
convergence expressions become more comprehensible. The result
requires a condition on the delay bound $\tau$ in terms of $n$ and the
ratio $\Lmax/\Lres$.
\begin{corollary} \label{co:thm_2}
Suppose that $\tau \ge 1$ and $n\ge 5$ and that
\beq \label{eq:boundtauc}
\tau(\tau + 1) \leq \frac{\sqrt{n} \Lmax}{4e\Lres}.
\eeq
If we choose
\beq \label{eq:choicerhoc}
\rho=1+{4e\tau \Lres\over  \sqrt{n}\Lmax},
\eeq
then the steplength $\gamma=1/2$ will satisfy the bounds
\eqnok{eq:boundgammac}. In addition, for the essentially strongly
convex case \eqnok{eq:esc} with $l>0$, we have for $j=1,2,\dotsc$ that
\begin{equation}
 \label{eqn_thm_2_good_c}
\E  (f(x_j)-f^*) 
 \le \left( 1-\frac{l}{n(l+2\Lmax)} \right)^j (\Lmax R_0^2 + f(x_0)-f^*),
\end{equation}
while for the case of general convex $f$, we have
\beq \label{eqn_thm_3_good_c}
\E (f(x_j)-f^*) \le  \frac{n(\Lmax R_0^2 + f(x_0)-f^*)}{j+n}.
\eeq
\end{corollary}

Similarly to Section~\ref{sec_unconstrained}, and provided $\tau$
satisfies \eqnok{eq:boundtauc}, the convergence rate is not affected
appreciably by the delay bound $\tau$, and near-linear speedup can
be expected for multicore implementations when \eqnok{eq:boundtauc}
holds. This condition is more restrictive than \eqnok{eq:boundtau} in
the unconstrained case, but still holds in many problems for
interesting values of $\tau$. When $\Lres/\Lmax$ is bounded
independently of dimension, the maximal number of cores allowed is of
the the order of $n^{1/4}$, which is smaller than the
$O(n^{1/2})$ value obtained for the unconstrained case.

%

We conclude this section with another high-probability bound, whose
proof tracks that of Theorem~\ref{thm:co1}.
\begin{theorem} \label{thm:co2}
Suppose that the conditions of Corollary~\ref{co:thm_2} hold,
including the definitions of $\rho$ and $\psi$. Then for $\epsilon>0$
and $\eta\in(0,1)$, we have that 
\[
\mathbb{P}\left(f(x_j)-f^*\le \epsilon \right) \ge 1-\eta,
\]
provided that one of the following conditions holds: In the
essentially strongly convex case~\eqnok{eq:esc} with $l>0$, we require
\beq \label{eq:def2_j1}
j \geq \frac{n(l+2\Lmax)}{l} \left|{\log{\Lmax R_0^2 + f(x_0)-f^*\over \epsilon\eta}}\right|, \nonumber
\eeq
while in the general convex case, it suffices that
\beq \label{eq:def2_j2}
j \geq \frac{n(\Lmax R_0^2+f(x_0)-f^*)}{\epsilon\eta} - n. \nonumber
\eeq
\end{theorem}

\section{Experiments}\label{sec_exp}

We illustrate the behavior of two variants of the stochastic
coordinate descent approach on test problems constructed from several
data sets. Our interests are in the efficiency of multicore
implementations (by comparison with a single-threaded implementation)
and in performance relative to alternative solvers for the same
problems.

All our test problems have the form \eqnok{eqn_mainproblem}, with
either $\Omega = \mathbb{R}^n$ or $\Omega$ separable as in
\eqnok{eq:sep}.  The objective $f$ is quadratic, that is,
\[
f(x) = \frac12 x^TQx + c^Tx,
\]
with $Q$ symmetric positive definite.

Our implementation of $\ASCD$ is called DIMM-WITTED (or DW for
short). It runs on various numbers of threads, from 1 to 40, each
thread assigned to a single core in our 40-core Intel Xeon
architecture. Cores on the Xeon architecture are arranged into four
sockets --- ten cores per socket, with each socket having its own
memory. Non-uniform memory access (NUMA) means that memory accesses to
local memory (on the same socket as the core) are less expensive than
accesses to memory on another socket.  In our DW implementation, we
assign each socket an equal-sized ``slice'' of $Q$, a row
submatrix. The components of $x$ are partitioned between cores, each
core being responsible for updating its own partition of $x$ (though
it can read the components of $x$ from other cores). The components of
$x$ assigned to the cores correspond to the rows of $Q$ assigned to
that core's socket. Computation is grouped into ``epochs,''
where an epoch is defined to be the period of computation during which
each component of $x$ is updated exactly once. We use the parameter
$p$ to denote the number of epochs that are executed between
reordering (shuffling) of the coordinates of $x$. We investigate both
shuffling after every epoch ($p=1$) and after every tenth epoch
($p=10$).  Access to $x$ is lock-free, and updates are performed
asynchronously. This update scheme does not implement exactly the
``sampling with replacement'' scheme analyzed in previous sections,
but can be viewed as a high performance, practical adaptation of the
$\ASCD$ method.

To do each coordinate descent update, a thread must read the latest
value of $x$. Most components are already in the cache for that core,
so that it only needs to fetch those components recently changed. When
a thread writes to $x_i$, the hardware ensures that this $x_i$ is
simultaneously removed from other cores, signaling that they must
fetch the updated version before proceeding with their respective
computations.

Although DW is not a precise implementation of $\ASCD$, it largely
achieves the consistent-read condition that is assumed by the
analysis. Inconsistent read happens on a core only if the following
three conditions are satisfied simultaneously:
\begin{itemize}
\item A core does not finish reading recently changed coordinates of
  $x$ (note that it needs to read no more than $\tau$ coordinates);
\item Among these recently changed coordinates, modifications take
  place both to coordinates that {\em have been read} and that are
  {\em still to be read} by this core;
\item Modification of the already-read coordinates happens earlier
  than the modification of the still-unread coordinates.
\end{itemize}
Inconsistent read will occur only if at least two coordinates of $x$
are modified twice during a stretch of approximately $\tau$ updates to
$x$ (that is, iterations of Algorithm~\ref{alg_ascd}). For the DW
implementation, inconsistent read would require repeated updating of a
particular component in a stretch of approximately $\tau$ iterations
that straddles two epochs. This event would be rare, for typical
values of $n$ and $\tau$. Of course, one can avoid the inconsistent
read issue altogether by changing the shuffling rule slightly,
enforcing the requirement that no coordinate can be modified twice in
a span of $\tau$ iterations. From the practical perspective, this
change does not improve performance, and detracts from the simplicity
of the approach. From the theoretical perspective, however, the
analysis for the inconsistent-read model would be interesting and
meaningful, and we plan to study this topic in future work.


The first test problem {\tt QP} is an unconstrained, regularized least
squares problem constructed with synthetic data. It has the form
\begin{equation} \label{eq:uncqp}
\min_{x\in\mathbb{R}^n} \; f(x):={1\over 2}\|Ax-b\|^2 + {\alpha\over
  2}\|x\|^2.
\end{equation}
All elements of $A\in\R^{m\times n}$, the true model
$\tilde{x}\in\R^{n}$, and the observation noise vector
$\delta\in\R^{m}$ are generated in i.i.d. fashion from the Gaussian
distribution $\mathcal{N}(0,1)$, following which each column in $A$ is
scaled to have a Euclidean norm of $1$. The observation
$b\in\mathbb{R}^m$ is constructed from $A\tilde{x} + \delta
\|A\tilde{x}\|/(5 m)$.
We choose $m=6000$,
$n=20000$, and $\alpha=0.5$. We therefore have
$\Lmax=1+\alpha=1.5$ and
\[
\frac{\Lres}{\Lmax} \approx {1+\sqrt{n/m}+\alpha \over 1+\alpha} \approx 2.2.
\]
With this estimate, the
condition~\eqref{eq:boundtau} is satisfied when delay parameter $\tau$
is less than about $95$. In Algorithm~\ref{alg_ascd}, we set the
steplength parameter $\gamma$ to $1$, and we choose initial iterate to
be $x_0=\bold{0}$. We measure convergence of the residual norm $\|
\nabla f(x) \|$.

Our second problem {\tt QPc} is a bound-constrained version of
\eqnok{eq:uncqp}:
\begin{equation} \label{eq:bcqp}
\min_{x\in\mathbb{R}^n_+}\quad f(x):={1\over
  2}(x-\tilde{x})^T(A^TA+\alpha I)(x-\tilde{x}).
\end{equation}
The methodology for generating $A$ and $\tilde{x}$ and for choosing
the values of $m$, $n$, $\gamma$, and $x_0$ is the same as for
\eqnok{eq:uncqp}. We measure convergence via the residual $\|x-\PO(x-
\nabla f(x))\|$, where $\Omega$ is the nonnegative orthant
$\mathbb{R}^n_{+}$. At the solution of \eqnok{eq:bcqp}, about half the
components of $x$ are at their lower bound of $0$.

Our third and fourth problems are quadratic penalty functions
for linear programming relaxations of vertex cover problems on large
graphs. The vertex cover problem for an undirected graph with edge set $E$ and
vertex set $V$ can be written as a binary linear program:
\[
\min_{y \in \{0,1\}^{|V|}} \, \sum_{v \in V} y_v \quad \mbox{subject to} \;\;
y_u+y_v \ge  1, \quad \forall \, (u,v) \in E.
\]
By relaxing each binary constraint to the interval $[0,1]$,
introducing slack variables for the cover inequalities, we obtain a
problem of the form
\[
\min_{y_v \in [0,1], \, s_{uv} \in [0,1]} \, 
\sum_{v \in V} y_v \quad \mbox{subject to} \;\;
y_u+y_v -s_{uv} = 0, \quad \forall \, (u,v) \in E.
\]
This has the form
\[
\min_{x \in [0,1]^n} \, c^Tx \quad \mbox{subject to} \;\; Ax=b,
\]
for $n=|V|+|E|$. The test problem \eqnok{eq:vc} is a regularized
quadratic penalty reformulation of this linear program for some
penalty parameter $\beta$:
%
\begin{equation} \label{eq:vc}
\min_{x\in [0,1]^n} \quad
c^Tx+{\beta \over 2} \|Ax-b\|^2 + {1\over 2\beta} \|x\|^2,
\end{equation}
with $\beta=5$. 
Two test data sets {\tt Amazon} and {\tt DBLP} have dimensions
$n=561050$ and $n=520891$, respectively.

\begin{figure}
  \centering
\includegraphics[scale=0.30]{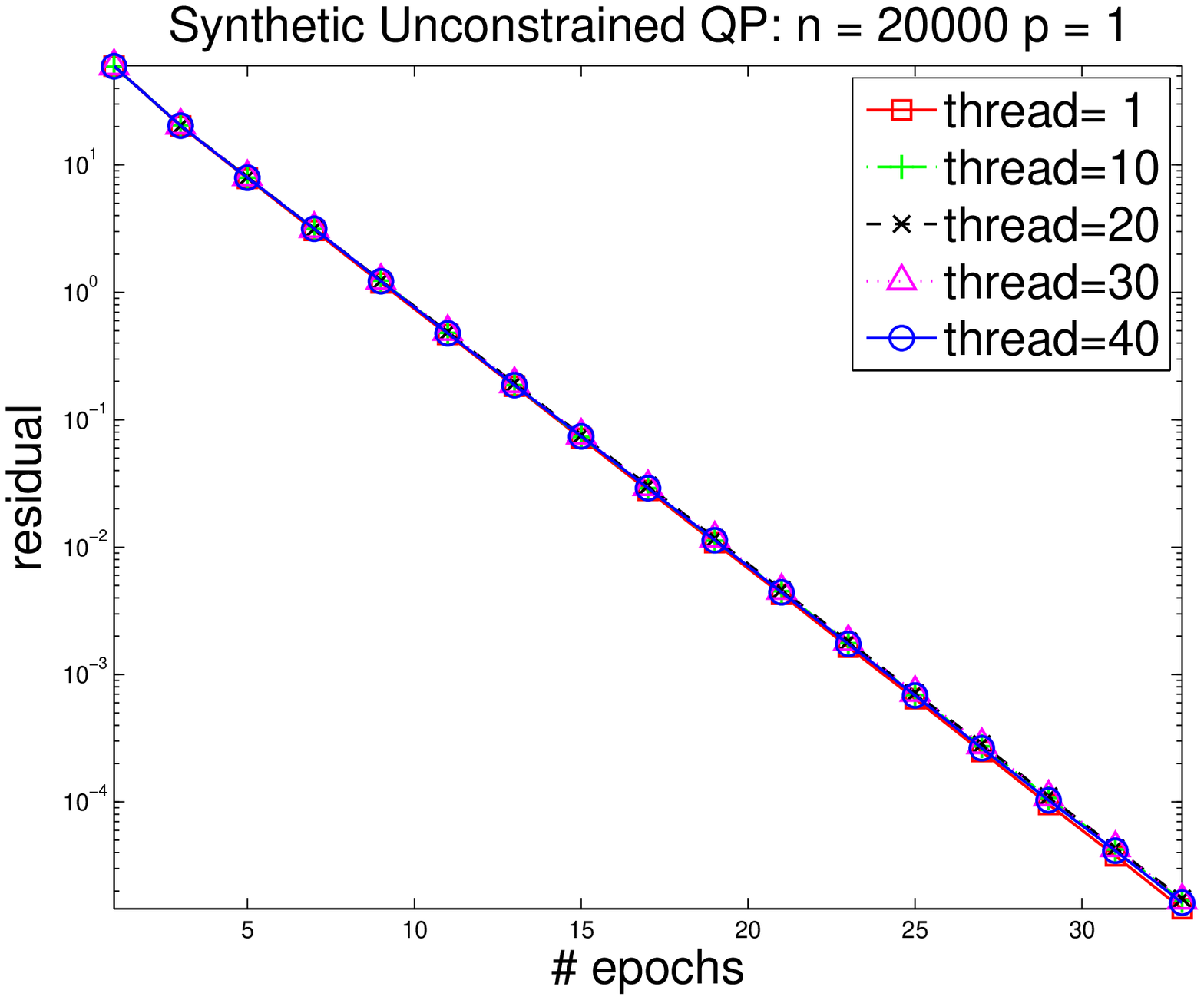}
\includegraphics[scale=0.30]{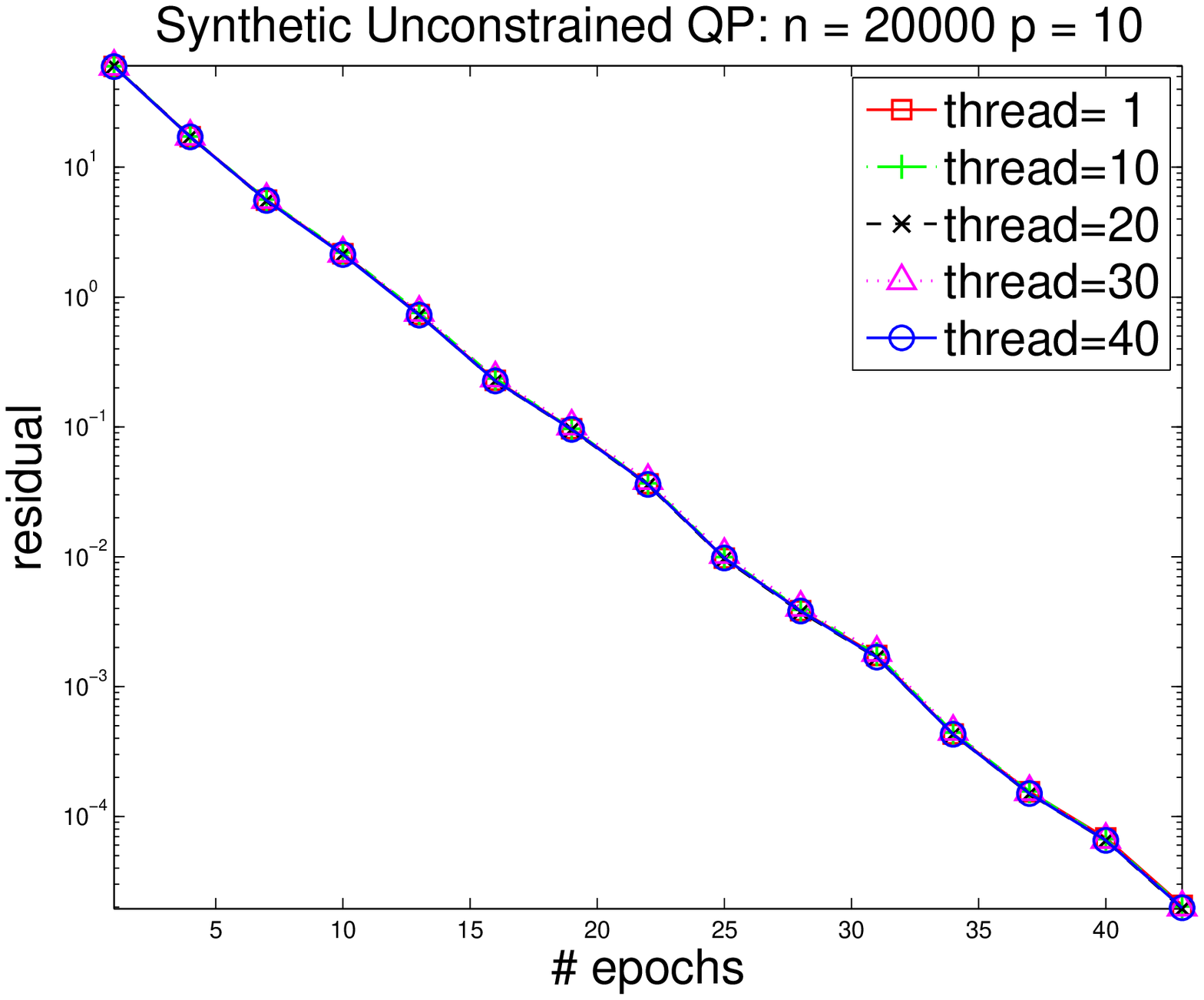}\\
\includegraphics[scale=0.30]{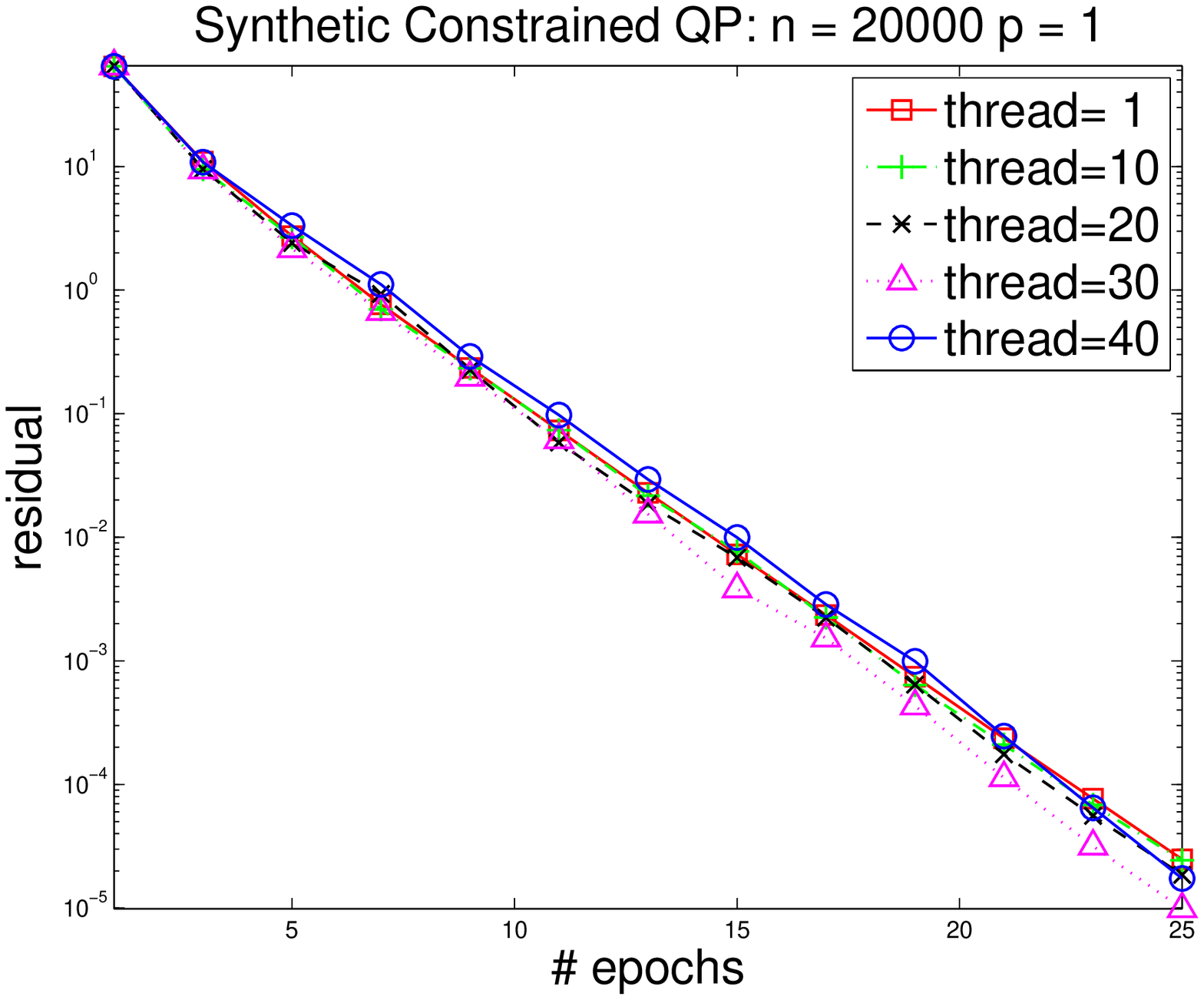}
\includegraphics[scale=0.30]{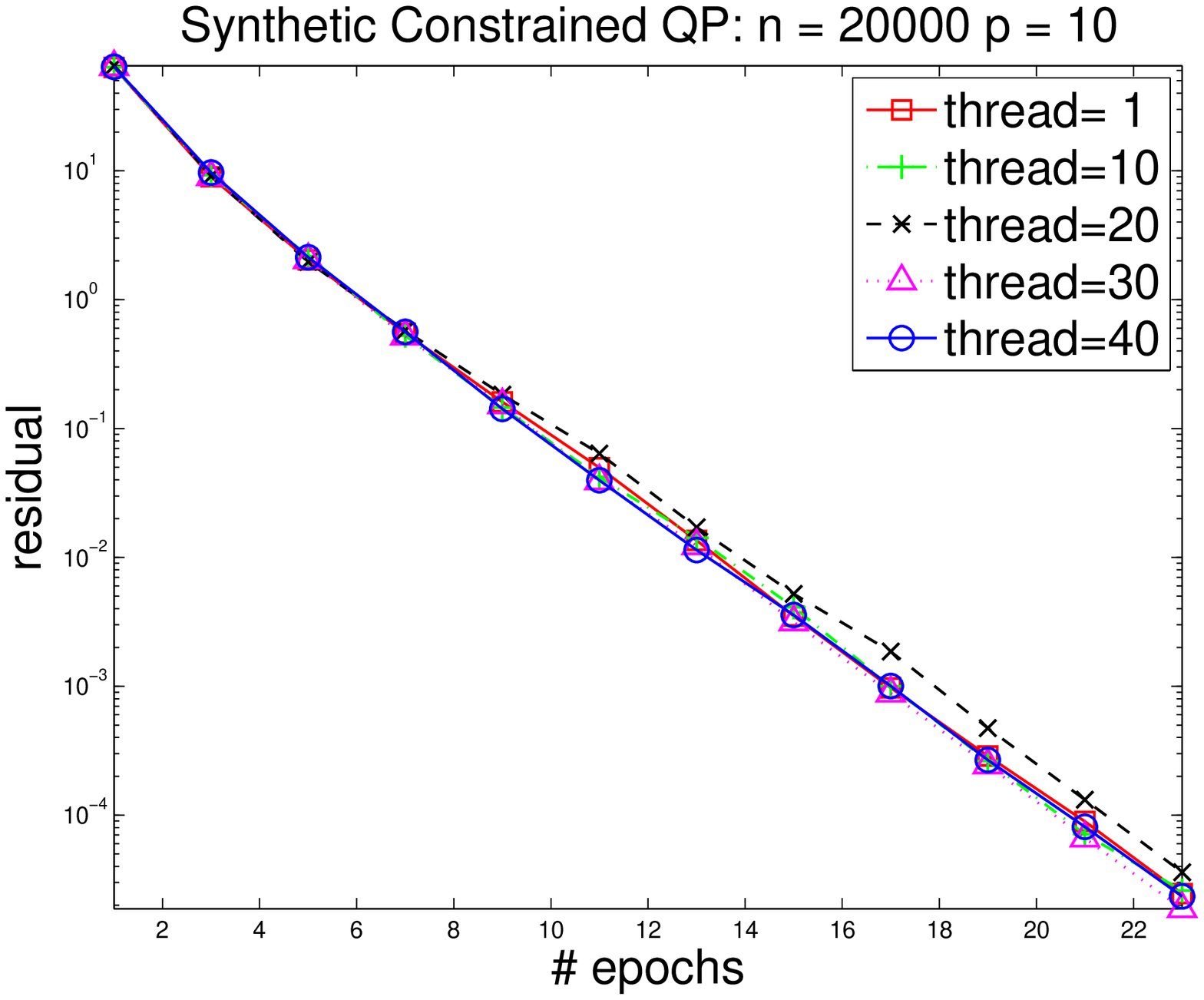}\\
\includegraphics[scale=0.30]{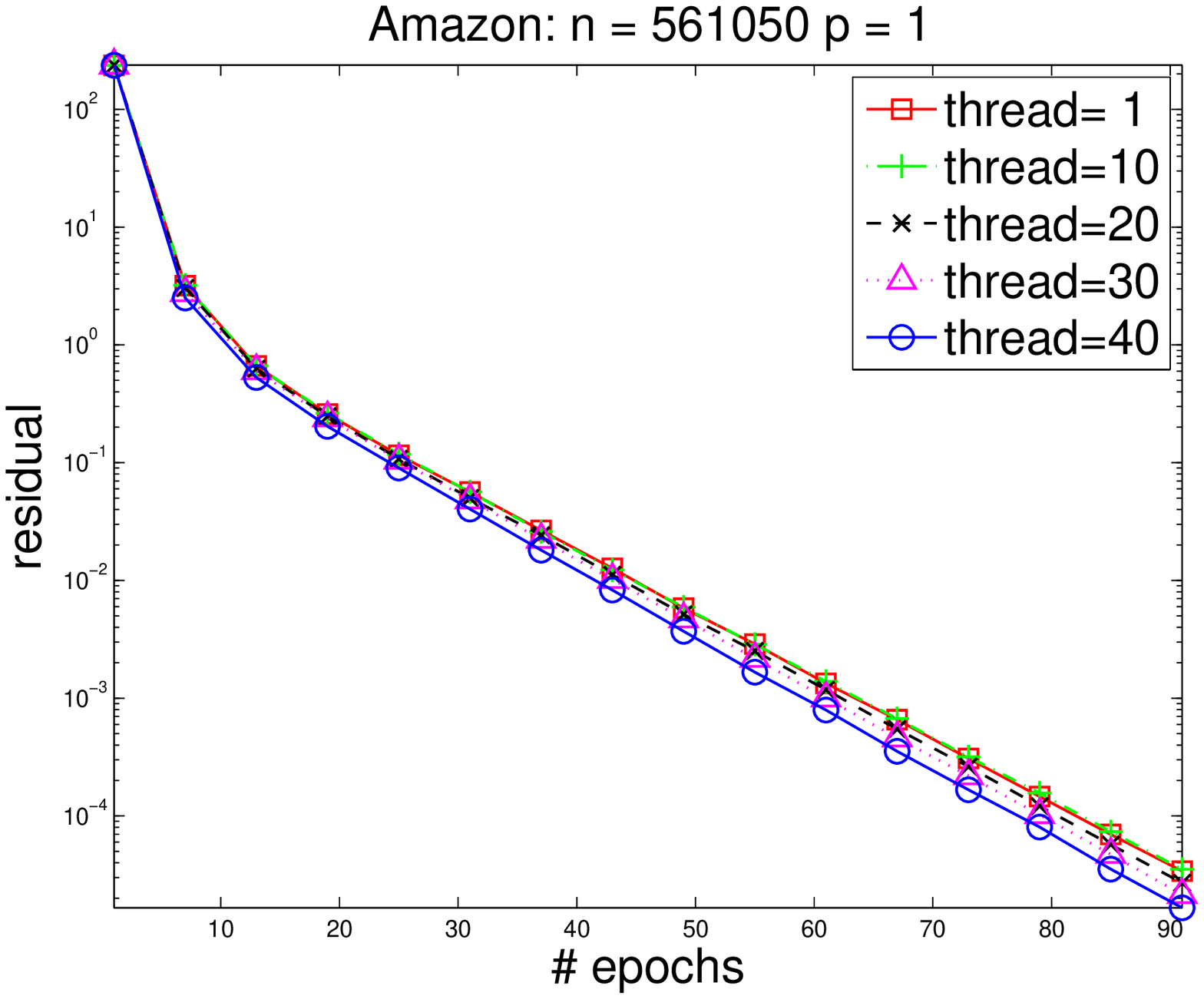}
\includegraphics[scale=0.30]{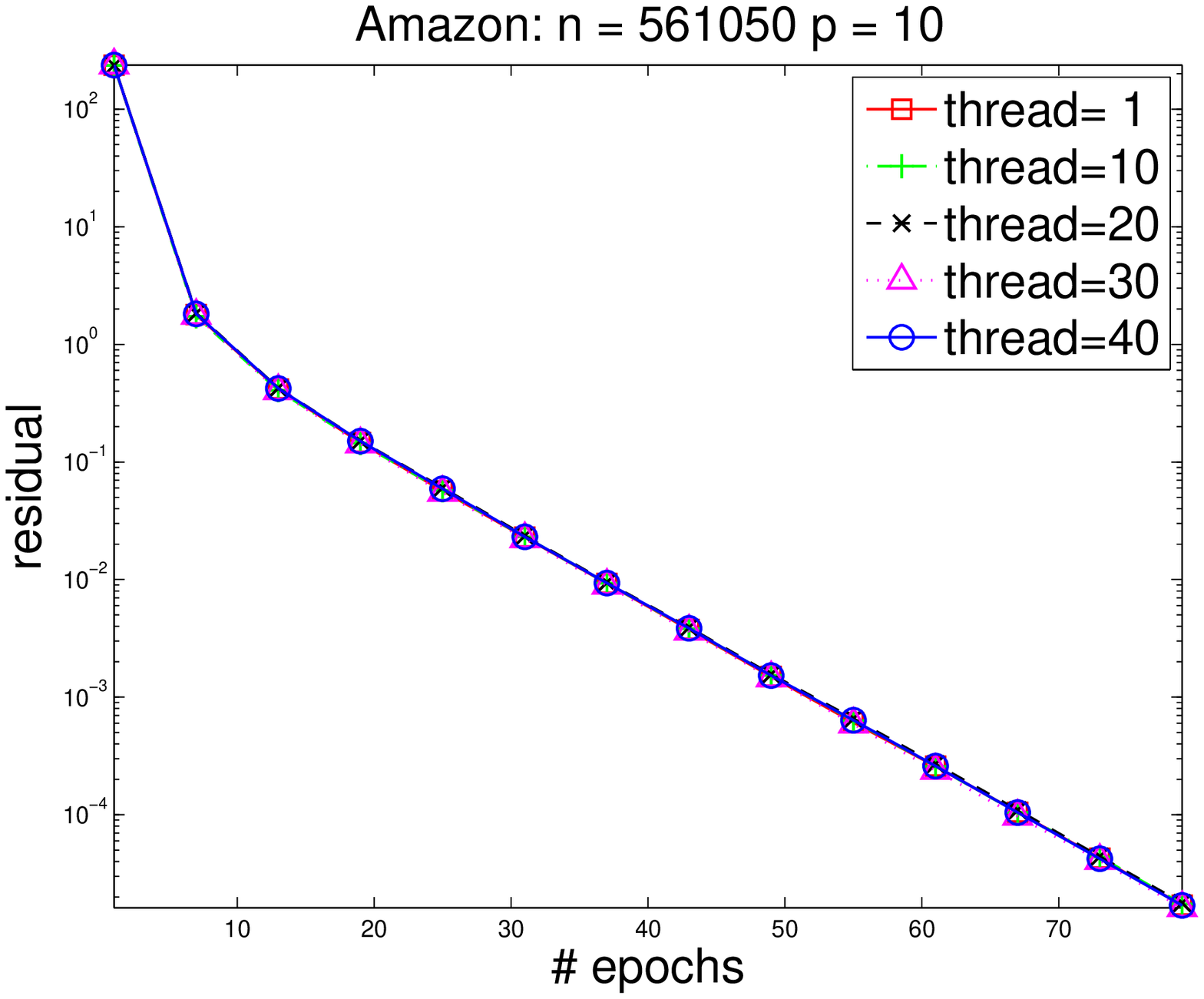}\\
\includegraphics[scale=0.30]{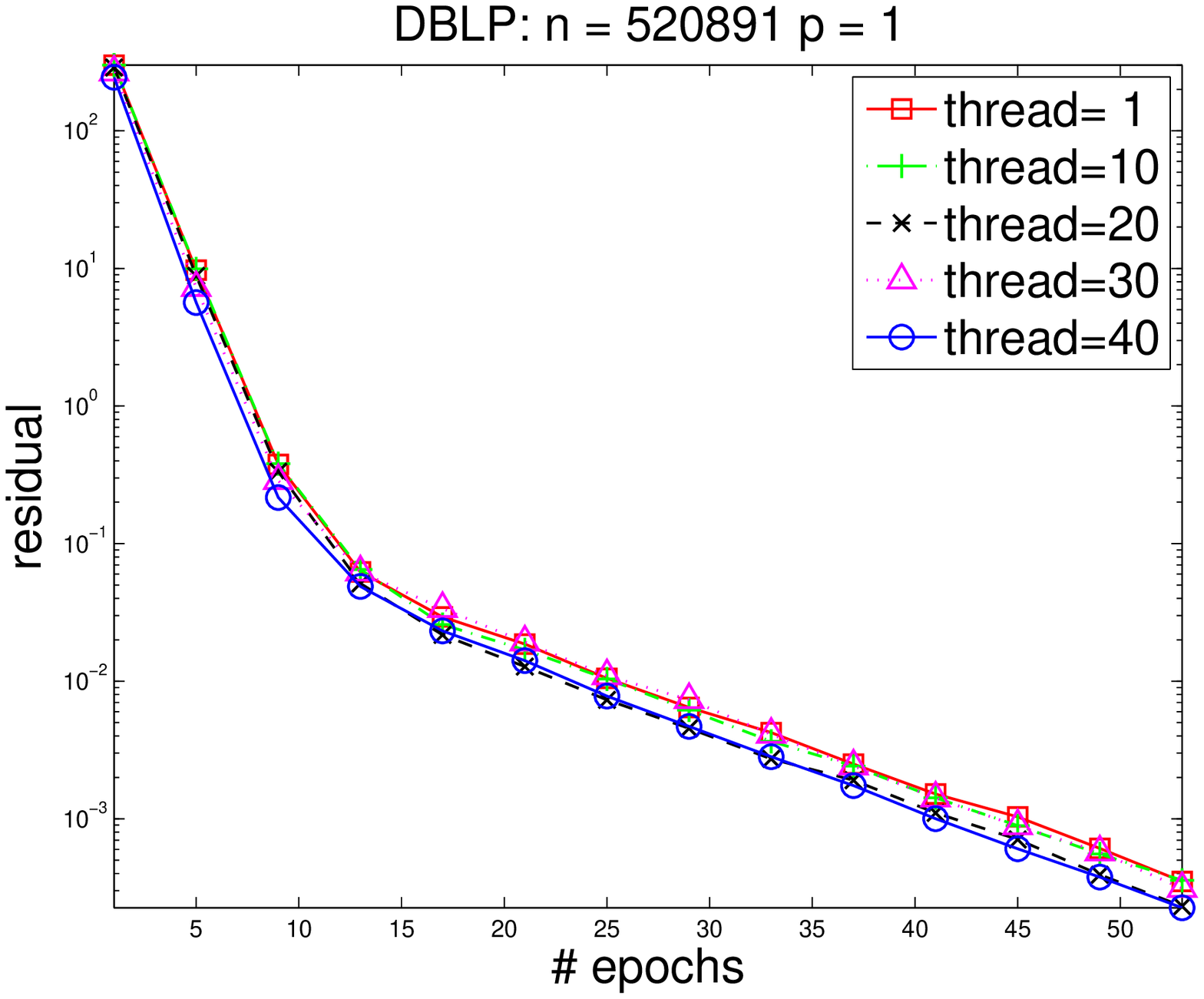}
\includegraphics[scale=0.30]{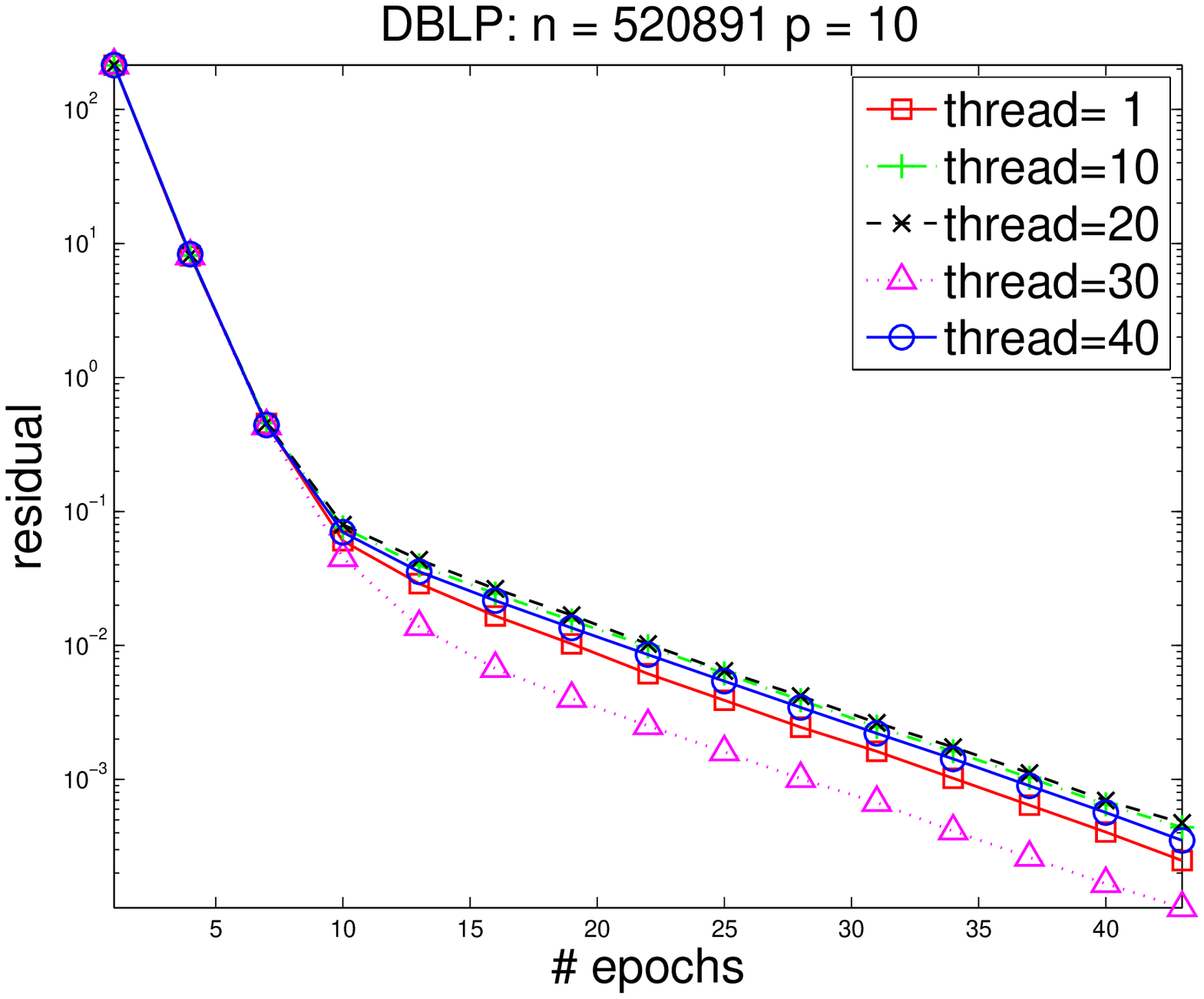}
\caption{Residuals vs epoch number for the four test problems. Results
  are reported for variants in which indices are reshuffled after
  every epoch ($p=1$) and after every tenth epoch ($p=10$).}
\label{fig:conv}
\end{figure}

We tracked the behavior of the residual as a function of the number of
epochs, when executed on different numbers of
cores. Figure~\ref{fig:conv} shows convergence behavior for each of
our four test problems on various numbers of cores with two different
shuffling periods: $p=1$ and $p=10$. We note the following points.
\begin{itemize}
\item The total amount of computation to achieve any level of
  precision appears to be almost independent of the number of cores,
  at least up to 40 cores. In this respect, the performance of the
  algorithm does not change appreciably as the number of cores is
  increased. Thus, any deviation from linear speedup is due not to
  degradation of convergence speed in the algorithm but rather to
  systems issues in the implementation. 
\item When we reshuffle after every epoch ($p=1$), convergence is
  slightly faster in synthetic unconstrained {\tt QP} but slightly
  slower in {\tt Amazon} and {\tt DBLP} than when we do occasional
  reshuffling ($p=10$). Overall, the convergence rates with different
  shuffling periods are comparable in the sense of epochs. However,
  when the dimension of the variable is large, the shuffling operation
  becomes expensive, so we would recommend using a large value for $p$
  for large-dimensional problems.
\end{itemize}

\begin{figure}[ht]
  \centering
 \includegraphics[scale=0.30]{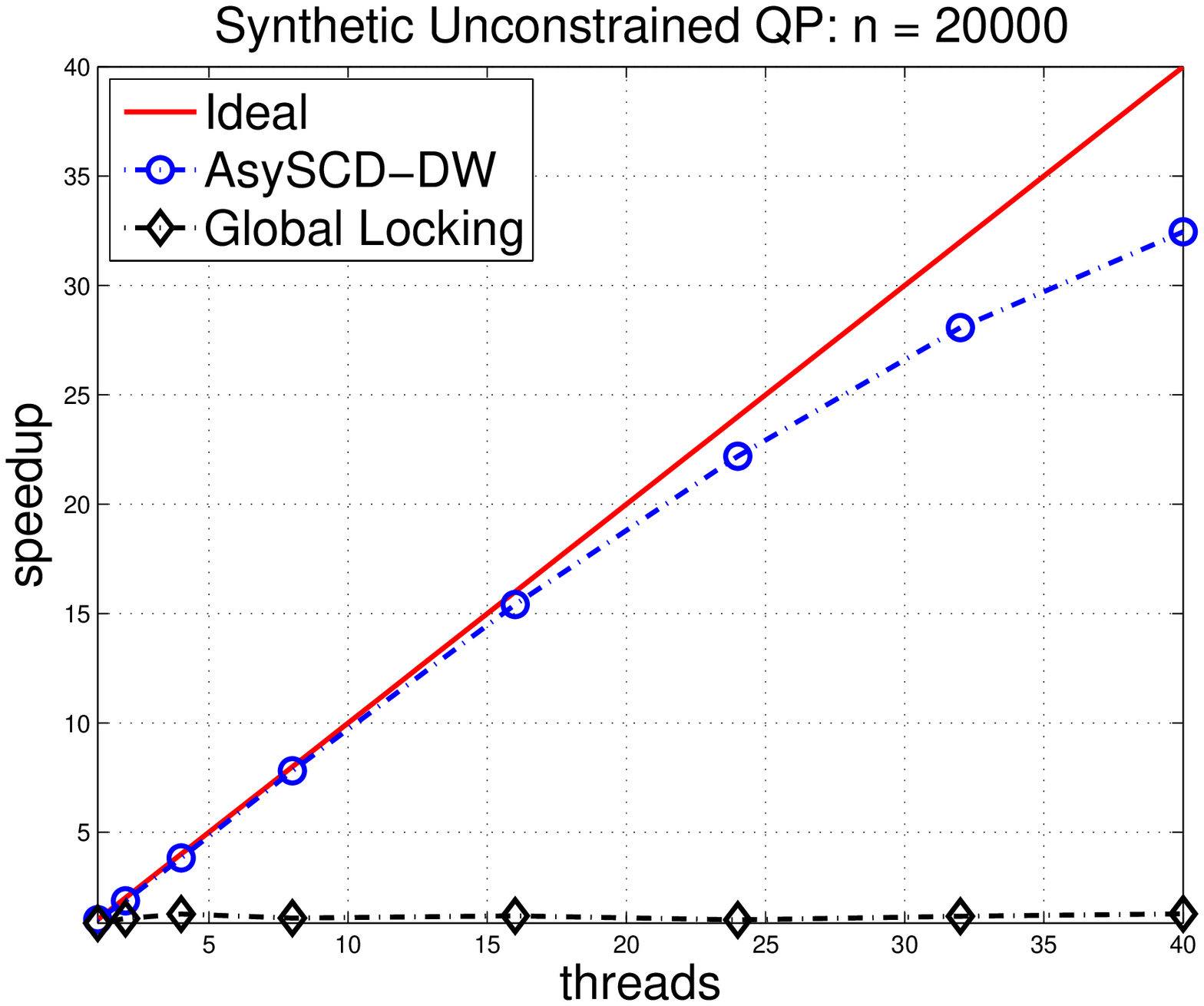}
\includegraphics[scale=0.30]{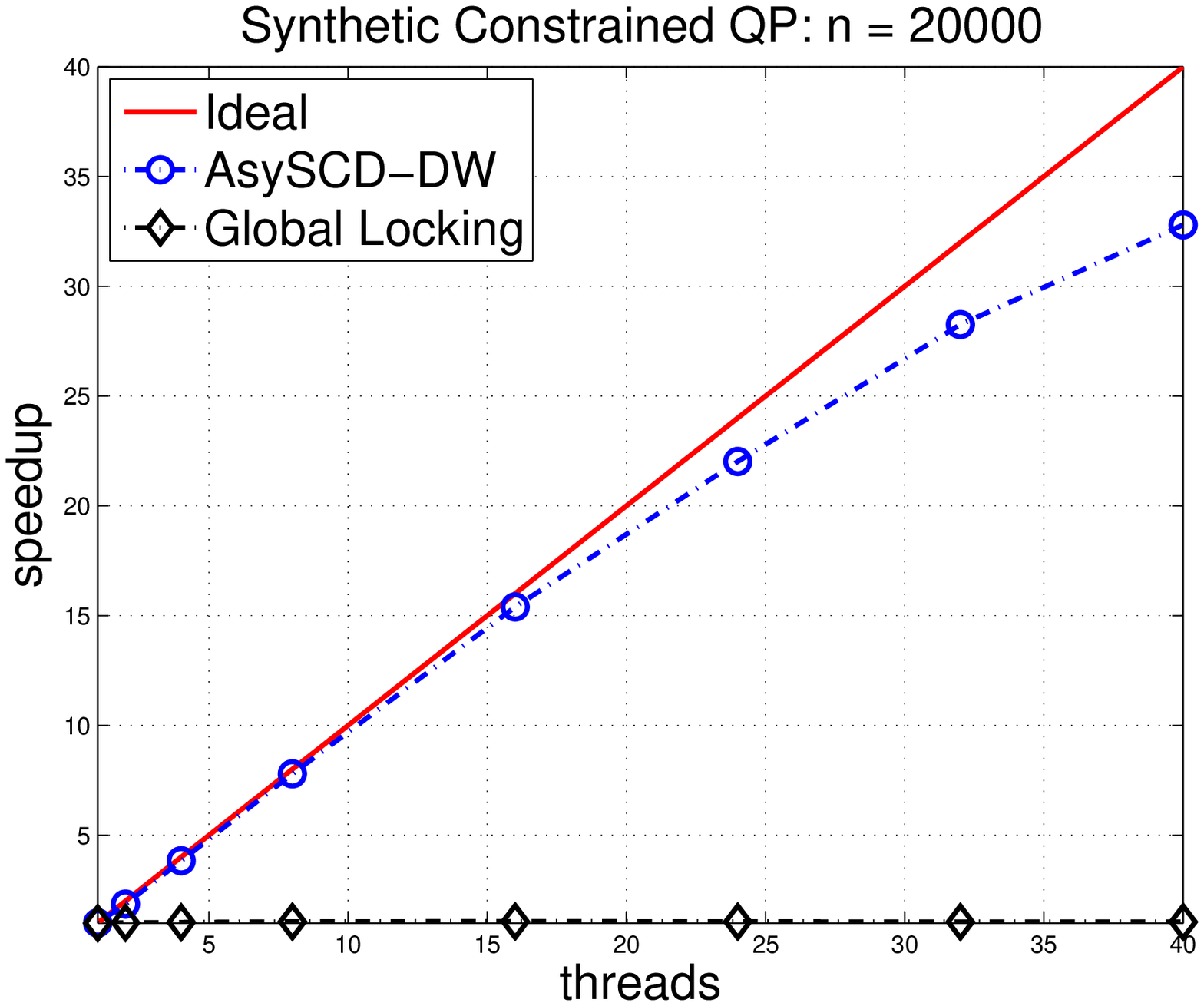}
       \caption{Test problems 1 and 2: Speedup of multicore
         implementations of DW on up to 40 cores of an Intel Xeon
         architecture. Ideal (linear) speedup curve is shown for
         reference, along with poor speedups obtained for a
         global-locking strategy.}
    \label{fig:speedup12}
\end{figure}

\begin{figure}[htp!]
  \centering
 \includegraphics[scale=0.30]{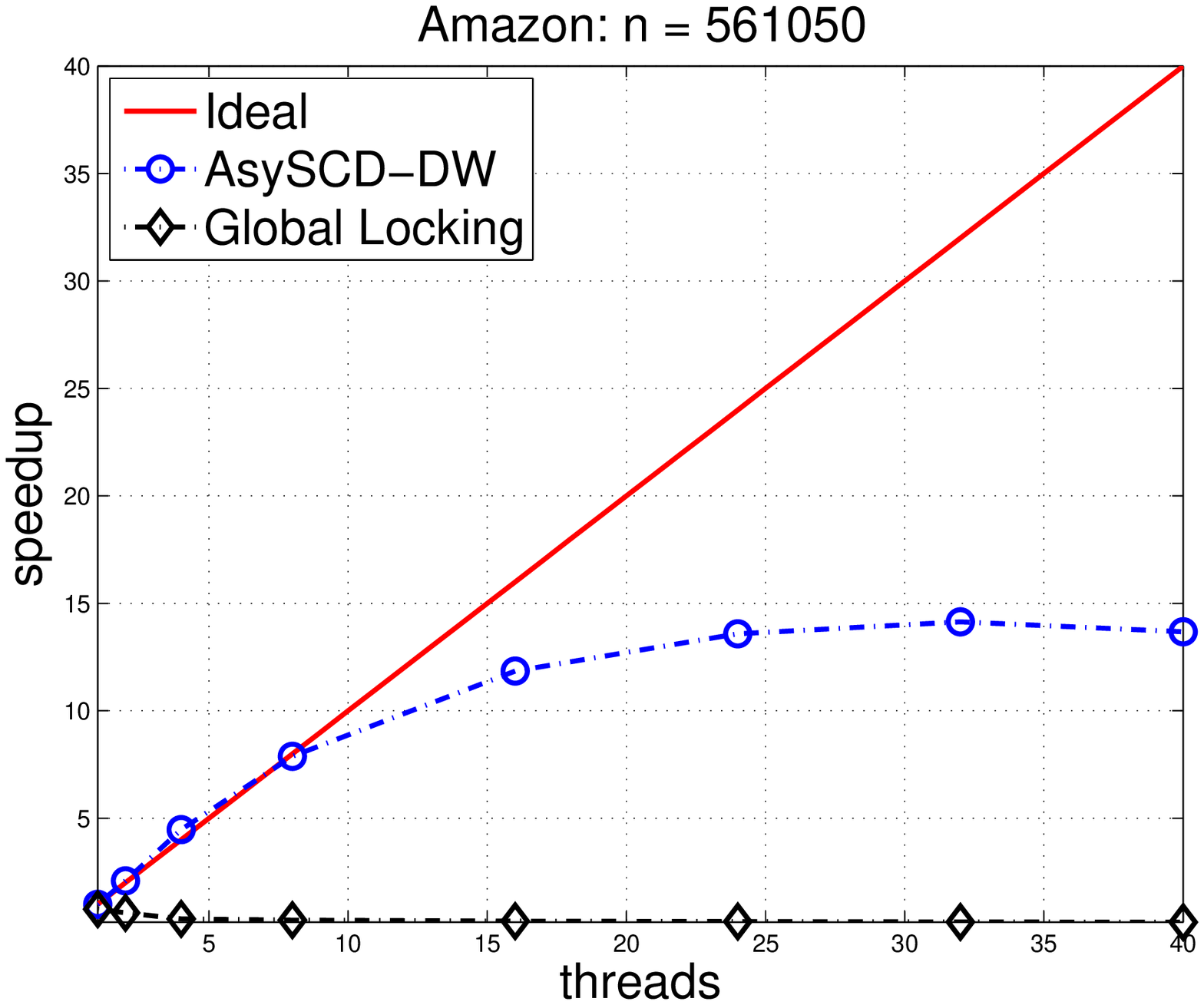}
\includegraphics[scale=0.30]{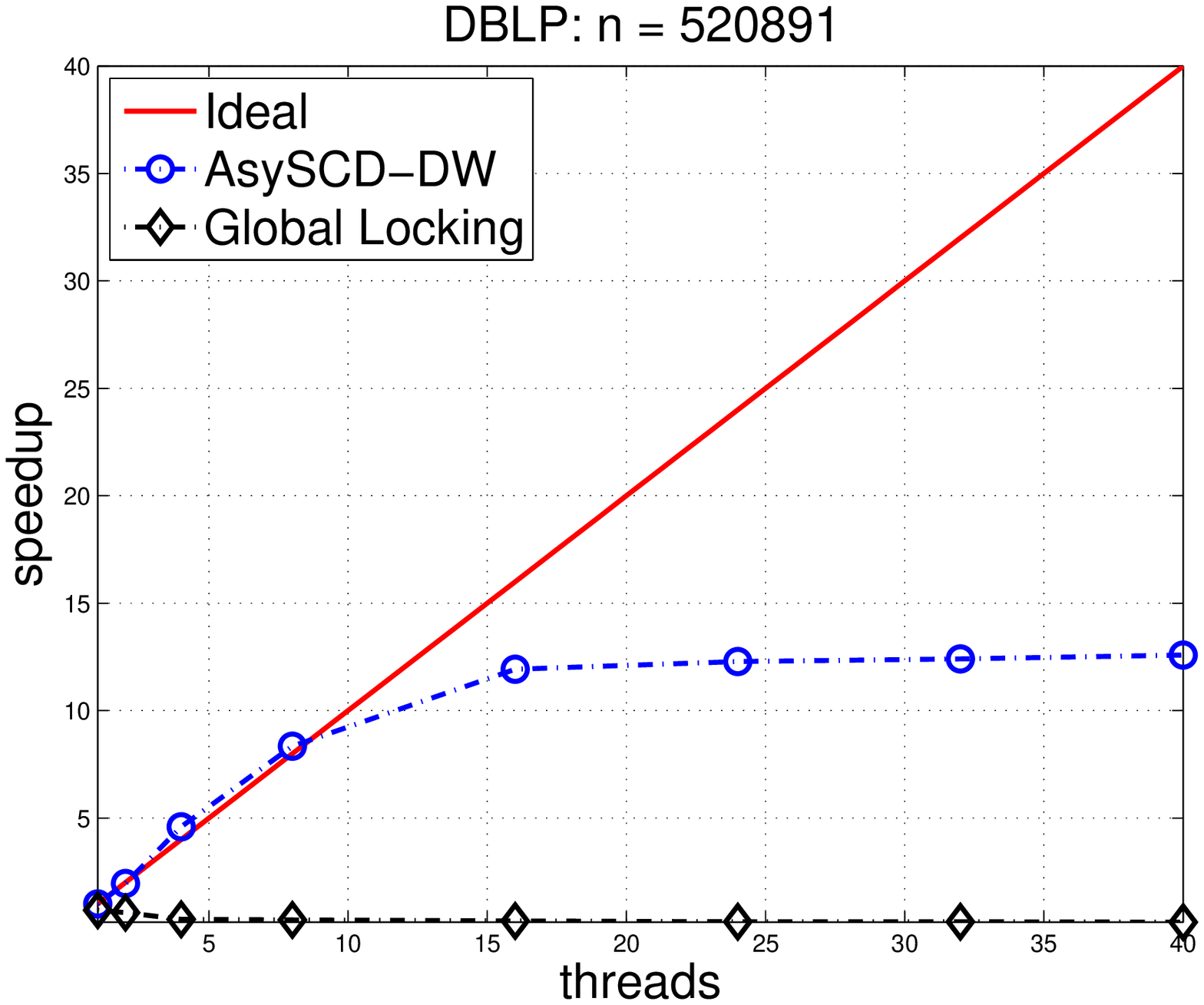}
       \caption{Test problems 3 and 4: Speedup of multicore implementations
  of DW on up to 40 cores of an Intel Xeon architecture. Ideal
  (linear) speedup curve is shown for reference, along with poor
  speedups obtained for a global-locking strategy.}
    \label{fig:speedup34}
\end{figure}

Results for speedup on multicore implementations are shown in
Figures~\ref{fig:speedup12} and \ref{fig:speedup34} for DW with
$p=10$.  Speedup is defined as follows:
\[
\frac{\text{runtime a single core using DW}}{\text{runtime on $P$ cores}}.
\]
Near-linear speedup can be observed for the two QP problems with
synthetic data. For Problems 3 and 4, speedup is at most 12-14; there
are few gains when the number of cores exceeds about 12. 
We believe that the degradation is due mostly to memory
contention. Although these problems have high dimension, the matrix
$Q$ is very sparse (in contrast to the dense $Q$ for the synthetic
data set). Thus, the ratio of computation to data movement / memory
access is much lower for these problems, making memory contention
effects more significant.

 Figures~\ref{fig:speedup12} and \ref{fig:speedup34} also show results
 of a global-locking strategy for the parallel stochastic coordinate
 descent method, in which the vector $x$ is locked by a core whenever
 it performs a read or update. The performance curve for this strategy
 hugs the horizontal axis; it is not competitive.



Wall clock times required for the four test problems on $1$ and $40$
cores, to reduce residuals below $10^{-5}$ are shown in
Table~\ref{ta:times}. (Similar speedups are noted when we use a
convergence tolerance looser than $10^{-5}$.)



\begin{table}[t!]
\centering
\begin{tabular}{|l|cc|}
\hline
Problem & 1 core & 40 cores \\ \hline
\texttt{QP} & 98.4 & 3.03 \\
\texttt{QPc} & 59.7 & 1.82 \\
\texttt{Amazon} & 17.1 & 1.25 \\
\texttt{DBLP} & 11.5 & .91 \\ \hline
\end{tabular}
\caption{Runtimes (seconds) for the four test problems on $1$ and $40$ cores.}
\label{ta:times}
\end{table}

\begin{table}[t!]
\centering
\begin{tabular}{|l|cc|}
\hline
\#cores 	& Time(sec) 	& Speedup	\\ \hline
1      	&   55.9     	& 1		\\ 
10      &   5.19        & 10.8   	\\ 
20      &   2.77        & 20.2   	\\ 
30      &   2.06        & 27.2   	\\ 
40      &   1.81        & 30.9	\\ \hline
\end{tabular}
\caption{Runtimes (seconds) and speedup for multicore implementations
  of DW on different number of cores for the weakly convex
  \texttt{QPc} problem (with $\alpha=0$) to achieve a residual below
  $0.06$.}
\label{ta:weaklyconvex}
\end{table}

\begin{table}[tp!]
\centering
\begin{tabular}{|l|cc|}
\hline
 \#cores 	&	Time(sec)  	&	Speedup	\\
 	& $\SynGD$ / $\ASCD$	& $\SynGD$ / $\ASCD$ 	\\
\hline
 1      &      96.8 / 27.1     &      0.28 / 1.00	\\
 10     &     11.4 / 2.57      &     2.38 / 10.5	\\
 20     &     6.00 / 1.36      &     4.51 / 19.9	\\
 30     &     4.44 /  1.01     &     6.10  / 26.8	\\	
 40     &     3.91 / 0.88      &     6.93 / 30.8	\\
\hline
\end{tabular}
\caption{Efficiency comparison between $\SynGD$ and $\ASCD$ for the
  {\tt QP} problem. The running time and speedup are based on the
  residual achieving a tolerance of $10^{-5}$.}
\label{ta:SynGD}
\end{table}

\begin{table} [tp!]
\centering
\begin{tabular}{ | c | r r | r r |} 
\hline
Dataset & \multicolumn{1}{c}{\# of} & \multicolumn{1}{c|}{\# of}
& \multicolumn{2}{c|}{Train time(sec)} \\ 
        & Samples & Features & LIBSVM & $\ASCD$ \\ 
\hline 
\textsf{adult} & 32561 & 123 & 16.15 & 1.39\\
\textsf{news} & 19996 & 1355191 & 214.48 & 7.22\\
\textsf{rcv} & 20242  & 47236 & 40.33 & 16.06\\
\textsf{reuters} & 8293 & 18930 & 1.63 & 0.81\\
\textsf{w8a} & 49749 & 300 & 33.62 & 5.86\\
\hline
\end{tabular}
\caption{Efficiency comparison between LIBSVM and $\ASCD$ for kernel
  SVM using 40 cores using homogeneous kernels ($K(x_i, x_j) = (x_i^Tx_j)^2$). The running time and speedup are
  calculated based on the ``residual'' $10^{-3}$. Here, to make both
  algorithms comparable, the ``residual'' is defined by
  $\|x-\mathcal{P}_{\Omega}(x-\nabla
  f(x))\|_\infty$.\label{ta:kernelSVM}}
\end{table}

All problems reported on above are essentially strongly convex. Similar speedup properties can be obtained in the weakly convex case as well. We show speedups for the {\tt QPc} problem with $\alpha=0$. Table~\ref{ta:weaklyconvex} demonstrates similar speedup to the essentially strongly convex case shown in Figure~\ref{fig:speedup12}.

Turning now to comparisons between $\ASCD$ and alternative algorithms,
we start by considering the basic gradient descent method. We
implement gradient descent in a parallel, synchronous fashion,
distributing the gradient computation load on multiple cores and
updating the variable $x$ in parallel at each step. The resulting
implementation is called $\SynGD$. Table~\ref{ta:SynGD} reports
running time and speedup of both $\ASCD$ over $\SynGD$, showing a
clear advantage for $\ASCD$. A high price is paid for the
synchronization requirement in $\SynGD$.

Next we compare $\ASCD$ to LIBSVM \citep{LIBSVM11} a popular multi-thread parallel solver for kernel support vector machines (SVM). Both algorithms are run on 40 cores to solve the dual formulation of kernel SVM, without an intercept term. All data sets used in~\ref{ta:kernelSVM} 
except \textsf{reuters} were obtained from the LIBSVM dataset repository\footnote{\url{http://www.csie.ntu.edu.tw/~cjlin/libsvmtools/datasets/}}. The dataset \textsf{reuters} is a sparse binary text classification dataset constructed as a one-versus-all version of Reuters-2159\footnote{\url{http://www.daviddlewis.com/resources/testcollections/reuters21578/}}. Our comparisons, shown in Table~\ref{ta:kernelSVM}, indicate that $\ASCD$ outperforms LIBSVM on these test sets.

\section{Extension}\label{sec_extension}

The $\ASCD$ algorithm can be extended by partitioning the coordinates
into blocks, and modifying Algorithm~\ref{alg_ascd} to work with these
blocks rather than with single coordinates.  If $L_i$, $\Lmax$, and
$\Lres$ are defined in the block sense, as follows:
\begin{align*}
\|\nabla f(x)-\nabla f(x+E_{i}t)\| &\leq \Lres\|t\|\quad\forall x, i, t\in\mathbb{R}^{|i|},\\
\|\nabla_i f(x)-\nabla_i f(x+E_{i}t)\| &\leq L_i\|t\|\quad \forall x, i, t\in\mathbb{R}^{|i|},\\
\Lmax &= \max_{i}L_i,
\end{align*}
where $E_i$ is the projection from the $i$th block to $\R^n$ and $|i|$
denotes the number of components in block $i$, our analysis can be
extended appropriately.

To make the $\ASCD$ algorithm more efficient, one can redefine the
steplength in Algorithm~\ref{alg_ascd} to be $\frac{\gamma}{L_{i(j)}}$
rather than $\frac{\gamma}{\Lmax}$.  Our analysis can be applied to
this variant by doing a change of variables to $\tilde{x}$, with
${x}_i=\frac{L_i}{\Lmax}\tilde{x}_i$ and defining $L_i$, $\Lres$,
and $\Lmax$ in terms of $\tilde{x}$.


\section{Conclusion}\label{sec_conclusion}

This paper proposes an asynchronous parallel stochastic coordinate
descent algorithm for minimizing convex objectives, in the
unconstrained and separable-constrained cases. Sublinear convergence
(at rate $1/K$) is proved for general convex functions, with stronger
linear convergence results for functions that satisfy an essential
strong convexity property. Our analysis indicates the extent to which
parallel implementations can be expected to yield near-linear speedup,
in terms of a parameter that quantifies the cross-coordinate
interactions in the gradient $\nabla f$ and a parameter $\tau$ that
bounds the delay in updating. Our computational experience confirms
the theory.


\acks{This project is supported by NSF Grants DMS-0914524,
  DMS-1216318, and CCF-1356918; NSF CAREER Award IIS-1353606; ONR
  Awards N00014-13-1-0129 and N00014-12-1-0041; AFOSR Award
  FA9550-13-1-0138; a Sloan Research Fellowship; and grants from
  Oracle, Google, and ExxonMobil.  }



\appendix
\section{Proofs for Unconstrained Case} \label{app:unc}

This section contains convergence proofs for $\ASCD$ in the
unconstrained case.



We start with a technical result, then move to the proofs of the three
main results of Section~\ref{sec_unconstrained}.

\begin{lemma} \label{lem_1}
For any $x$, we have
\[
{\|x-\PS(x)\|^2}\|\nabla f(x)\|^2 \geq {(f(x)-f^*)^2}.
\]
If the essential strong convexity property \eqnok{eq:esc} holds, we have
\[
\|\nabla f(x)\|^2 \geq 2l (f(x) - f^*).
\]
\end{lemma}
\begin{proof}
The first inequality is proved as follows:
\[
f(x)-f^* \leq \langle \nabla f(x), x-\PS(x) \rangle \leq \|\nabla f(x)\|\|\PS(x)-x\|.
\]
For the second bound, we have from the definition \eqnok{eq:esc},
setting $y \leftarrow x$ and $x \leftarrow \PS(x)$, that
\begin{align*}
f^* - f(x) & \geq \langle \nabla f(x), \PS(x)-x \rangle + {l\over 2} \|x-\PS(x)\|^2 \\
&=  {l \over 2} \|\PS(x)-x + \frac{1}{l} \nabla f(x)\|^2 - {1\over 2l}\|\nabla f(x)\|^2 \\
& \geq - {1\over 2l}\|\nabla f(x)\|^2,
\end{align*}
as required.
\end{proof}


\begin{proof} (Theorem~\ref{AsySCD:thm_1})
We prove each of the two inequalities in \eqref{eqn_thm_1} by
induction. We start with the left-hand inequality.  For all values of
$j$, we have
\begin{align}
\nonumber
&\E\left(\|\nabla f(x_j)\|^2 - \|\nabla f(x_{j+1})\|^2\right) \\
\nonumber
& \quad = \E \langle \nabla f(x_j)+\nabla f(x_{j+1}), \nabla f(x_j) - \nabla f(x_{j+1})\rangle \\
\nonumber
&\quad = \E \langle 2\nabla f(x_j)+\nabla f(x_{j+1}) - \nabla f(x_j), \nabla f(x_j) - \nabla f(x_{j+1})\rangle \\
\nonumber
&\quad \leq 2\E \langle \nabla f(x_j), \nabla f(x_j) - \nabla f(x_{j+1})\rangle \\
\nonumber
&\quad \leq 2\E (\|\nabla f(x_j)\|\|\nabla f(x_j) - \nabla f(x_{j+1})\|)\\
\nonumber
&\quad \leq 2\Lres\E (\|\nabla f(x_j)\|\|x_j - x_{j+1}\|) \\
\nonumber
&\quad \leq \frac{2\Lres\gamma}{\Lmax}\E (\|\nabla f(x_j)\|\|\nabla_{i(j)}f(x_{k(j)})\|) \\
\nonumber
&\quad \leq \frac{\Lres\gamma}{\Lmax}\E (n^{-1/2}\|\nabla f(x_j)\|^2+n^{1/2}\|\nabla_{i(j)}f(x_{k(j)})\|^2) \\
\nonumber
&\quad = \frac{\Lres\gamma}{\Lmax}\E (n^{-1/2}\|\nabla f(x_j)\|^2+n^{1/2}\E_{i(j)}(\|\nabla_{i(j)}f(x_{k(j)})\|^2)) \\
\nonumber
&\quad = \frac{\Lres\gamma}{\Lmax}\E (n^{-1/2}\|\nabla f(x_j)\|^2+n^{-1/2}\|\nabla f(x_{k(j)})\|^2) \\
\label{eqn_proof_thm_3}
&\quad \leq {\Lres\gamma \over \sqrt{n}{\Lmax}} \E \left(\|\nabla f(x_j)\|^2+ \|\nabla f(x_{k(j)})\|^2\right).
\end{align}
We can use this bound to show that the left-hand inequality in
\eqnok{eqn_thm_1} holds for $j=0$. By setting $j=0$ in
\eqnok{eqn_proof_thm_3} and noting that $k(0)=0$, we obtain
\beq \label{eq:goof1}
\E \left( \| \nabla f(x_0) \|^2 - \| \nabla f (x_1)\|^2 \right)
\le \frac{\Lres \gamma}{\sqrt{n} \Lmax} 2 \E ( \| \nabla f(x_0) \|^2).
\eeq
From 
\eqnok{eq:boundgamma.2}, 
we have
\[
\frac{2\Lres \gamma}{\sqrt{n} \Lmax} \le \frac{\rho-1}{\rho^{\tau}} \le
\frac{\rho-1}{\rho} = 1-\rho^{-1},
\]
where the second inequality follows from $\rho>1$. By substituting into
\eqnok{eq:goof1}, we obtain $\rho^{-1} \E (\| \nabla f(x_0) \|^2) \le
\E (\| \nabla f(x_1) \|^2)$, establishing the result for $j=1$. For
the inductive step, we use \eqnok{eqn_proof_thm_3} again, assuming
that the left-hand inequality in \eqnok{eqn_thm_1} holds up to stage
$j$, and thus that
\[
\E ( \| \nabla f(x_{k(j)}) \|^2) \le \rho^{\tau} \E ( \| \nabla f(x_j) \|^2),
\]
provided that $0 \le j-k(j) \le \tau$, as assumed. By substituting into
the right-hand side of \eqnok{eqn_proof_thm_3} again, and using
$\rho>1$, we obtain
\[
\E\left(\|\nabla f(x_j)\|^2 - \|\nabla f(x_{j+1})\|^2\right) \le
\frac{2\Lres \gamma \rho^{\tau}}{\sqrt{n} \Lmax} \E \left( \| \nabla
f(x_j) \|^2 \right).
\]
By substituting 
\eqnok{eq:boundgamma.2}
we conclude that the left-hand inequality in
\eqnok{eqn_thm_1} holds for all $j$.


We now work on the right-hand inequality in \eqnok{eqn_thm_1}.  For
all $j$, we have the following:
\begin{align}
\nonumber
&\E\left(\|\nabla f(x_{j+1})\|^2 - \|\nabla f(x_{j})\|^2\right) \\
\nonumber
& \quad = \E \langle \nabla f(x_j)+\nabla f(x_{j+1}), \nabla f(x_{j+1}) - \nabla f(x_{j})\rangle \\
\nonumber
&\quad \leq \E (\|\nabla f(x_j)+\nabla f(x_{j+1})\|\|\nabla f(x_j) - \nabla f(x_{j+1})\|)\\
\nonumber
&\quad \leq \Lres\E (\|\nabla f(x_j)+\nabla f(x_{j+1})\|\|x_j - x_{j+1}\|) \\
\nonumber
&\quad \leq \Lres\E ((2\|\nabla f(x_j)\|+\|\nabla f(x_{j+1})-\nabla f(x_j)\|)\|x_j - x_{j+1}\|) \\
\nonumber
&\quad \leq \Lres\E (2\|\nabla f(x_j)\| \|x_j - x_{j+1}\|+\Lres\|x_j - x_{j+1}\|^2) \\
\nonumber
&\quad \leq \Lres\E \left(\frac{2\gamma}{\Lmax}\|\nabla f(x_j)\| \|\nabla_{i(j)}f(x_{k(j)})\|+\frac{\Lres\gamma^2}{\Lmax^2}\|\nabla_{i(j)} f(x_{k(j)})\|^2\right) \\
\nonumber
&\quad \leq \Lres\E \bigg(\frac{\gamma}{\Lmax}(n^{-1/2}\|\nabla f(x_j)\|^2 + n^{1/2}\|\nabla_{i(j)}f(x_{k(j)})\|^2+\frac{\Lres\gamma^2}{\Lmax^2}\|\nabla_{i(j)} f(x_{k(j)})\|^2\bigg) \\
\nonumber
&\quad = \Lres\E \bigg(\frac{\gamma}{\Lmax}(n^{-1/2}\|\nabla f(x_j)\|^2 + n^{1/2}\E_{i(j)}(\|\nabla_{i(j)}f(x_{k(j)})\|^2))+\\
\nonumber
&\quad\quad \frac{\Lres\gamma^2}{\Lmax^2}\E_{i(j)}(\|\nabla_{i(j)} f(x_{k(j)})\|^2) \bigg) \\
\nonumber
&\quad = \Lres\E \left(\frac{\gamma}{\Lmax}(n^{-1/2}\|\nabla f(x_j)\|^2 + n^{-1/2}\|\nabla f(x_{k(j)})\|^2)+\frac{\Lres\gamma^2}{n\Lmax^2}\|\nabla f(x_{k(j)})\|^2\right) \\
\nonumber
&\quad = {\gamma \Lres\over \sqrt{n}{\Lmax}} \E \left(\|\nabla f(x_j)\|^2+ \|\nabla f(x_{k(j)})\|^2\right)+{\gamma^2\Lres^2\over n{\Lmax^2}}\E (\|\nabla f(x_{k(j)})\|^2) \\
\label{eqn_proof_thm_4}
&\quad \leq {\gamma \Lres\over \sqrt{n}{\Lmax}} \E (\|\nabla f(x_j)\|^2)+ \left({\gamma \Lres\over \sqrt{n}{\Lmax}}+{\gamma \Lres^2\over n{\Lmax^2}}\right) \E(\|\nabla f(x_{k(j)})\|^2),
\end{align}
where the last inequality is from the observation $\gamma\leq 1$.
By setting $j=0$ in this bound, and noting that $k(0)=0$, we obtain
\beq \label{eq:goof2}
\E\left(\|\nabla f(x_1)\|^2 - \|\nabla f(x_0)\|^2\right)  \le
\left( \frac{2\gamma \Lres}{\sqrt{n} \Lmax} + \frac{\gamma \Lres^2}{n \Lmax^2} \right)
\E (\| \nabla f(x_0) \|^2 ).
\eeq
By using 
 \eqnok{eq:boundgamma.3}, we have
\[
\frac{2\gamma \Lres}{\sqrt{n} \Lmax} + \frac{\gamma \Lres^2}{n \Lmax^2}
= \frac{\Lres \gamma}{\sqrt{n} \Lmax}
\left( 2 + \frac{\Lres}{\sqrt{n} \Lmax} \right)
\le \frac{\rho-1}{\rho^{\tau}} < \rho-1,
\]
where the last inequality follows from $\rho>1$. By substituting into
\eqnok{eq:goof2}, we obtain $\E ( \| \nabla f(x_1) \|^2) \le \rho \E (
\| \nabla f(x_0) \|^2)$, so the right-hand bound in \eqnok{eqn_thm_1}
is established for $j=0$.  For the inductive step, we use
\eqnok{eqn_proof_thm_4} again, assuming that the right-hand inequality
in \eqnok{eqn_thm_1} holds up to stage $j$, and thus that
\[
\E ( \| \nabla f(x_j) \|^2) \le \rho^{\tau} \E ( \| \nabla f(x_{k(j)}) \|^2),
\]
provided that $0 \le j-k(j) \le \tau$, as assumed. From
\eqnok{eqn_proof_thm_4} and the left-hand inequality in~\eqref{eqn_thm_1}, we have by substituting this bound that
\beq \label{eq:goof3}
\E\left(\|\nabla f(x_{j+1})\|^2 - \|\nabla f(x_{j})\|^2\right) \le
 \left( \frac{2\gamma \Lres\rho^{\tau}}{\sqrt{n}{\Lmax}} +
\frac{\gamma \Lres^2 \rho^{\tau}}{n \Lmax^2} \right) \E (\|\nabla
f(x_j)\|^2).
\eeq
It follows immediately from 
\eqnok{eq:boundgamma.3} 
that the term in parentheses in \eqnok{eq:goof3}
is bounded above by $\rho-1$. By substituting this bound into
\eqnok{eq:goof3}, we obtain $\E ( \| \nabla f(x_{j+1}) \|^2) \le \rho
\E ( \| \nabla f(x_j) \|^2)$, as required.

At this point, we have shown that both inequalities in
\eqnok{eqn_thm_1} are satisfied for all $j$.

Next we prove \eqref{eqn_thm_2} and \eqref{eqn_thm_3}. Take the
expectation of $f(x_{j+1})$ in terms of $i(j)$:
\begin{align}
\nonumber
&~~\E_{i(j)} f(x_{j+1}) \\
\nonumber
& = \E_{i(j)} f \left( x_j-\frac{\gamma}{\Lmax} e_{i(j)} \nabla_{i(j)}f(x_{k(j)}) \right)\\
\nonumber
& = {1\over n}\sum_{i=1}^n f \left( x_j-\frac{\gamma}{\Lmax} e_i\nabla_i f(x_{k(j)}) \right) \\
\nonumber
& \leq {1\over n} \sum_{i=1}^n f(x_j) - \frac{\gamma}{\Lmax}\langle \nabla f(x_j), e_i\nabla_i f(x_{k(j)})\rangle + {L_i\over 2 \Lmax^2}\gamma^2 \|\nabla_i f(x_{k(j)})\|^2 \\
\nonumber
& \leq f(x_j) - {\gamma\over n{\Lmax}}\langle \nabla f(x_j), \nabla f(x_{k(j)}) \rangle + {\gamma^2\over 2n{\Lmax}} \|\nabla f(x_{k(j)})\|^2
\\
\nonumber
& = f(x_j) + {\gamma\over n{\Lmax}} \underbrace{\langle \nabla f(x_{k(j)})-\nabla f(x_j), \nabla f(x_{k(j)})\rangle}_{T_1} \\
&\quad\quad - \left({\gamma\over n{\Lmax}} - {\gamma^2\over 2n{\Lmax}}\right)\|\nabla f(x_{k(j)})\|^2. \label{eqn_proof_thm_1}
\end{align}

The second term $T_1$ is caused by delay. If there is no delay, $T_1$ should be $0$ because of $\nabla f(x_{j}) = \nabla
f(x_{k(j)})$. We estimate the upper bound of $\|\nabla
f(x_{k(j)})-\nabla f(x_{j})\|$:
\begin{align}
\nonumber
\|\nabla f(x_{k(j)})-\nabla f(x_{j})\| & \leq \sum_{d=k(j)}^{j-1}\left\| \nabla f(x_{d+1})- \nabla f(x_d)\right\|\\
\nonumber
& \leq \Lres\sum_{d=k(j)}^{j-1} \left\|x_{d+1}-x_d\right\| \\
\label{eqn_proof_thm_1_2}
& = \frac{\Lres\gamma}{\Lmax} \sum_{d=k(j)}^{j-1} \left\|\nabla_{i(d)} f(x_{k(d)})\right\|.
\end{align}
Then $\E  (|T_1|)$ can be bounded by
\begin{align}
\nonumber
\E (|T_1|) &\leq \E (\|\nabla f(x_{k(j)})-\nabla f(x_{j})\|\|\nabla f(x_{k(j)})\|) \\
\nonumber
& \leq \frac{\Lres\gamma}{\Lmax} \E \left(\sum_{d=k(j)}^{j-1} \|\nabla_{i(d)} f(x_{k(d)})\| \|\nabla f(x_{k(j)})\|\right) \\
\nonumber
& \leq \frac{\Lres\gamma}{2\Lmax} \E \left(\sum_{d=k(j)}^{j-1} n^{1/2}\|\nabla_{i(d)} f(x_{k(d)})\|^2 + n^{-1/2} \|\nabla f(x_{k(j)})\|^2\right) \\
\nonumber
& = \frac{\Lres\gamma}{2\Lmax} \E \left(\sum_{d=k(j)}^{j-1} n^{1/2}\E_{i(d)}(\|\nabla_{i(d)} f(x_{k(d)})\|^2) + n^{-1/2} \|\nabla f(x_{k(j)})\|^2\right) \\
\nonumber
& = \frac{\Lres\gamma}{2\Lmax} \E \left(\sum_{d=k(j)}^{j-1} n^{-1/2}\|\nabla f(x_{k(d)})\|^2 + n^{-1/2} \|\nabla f(x_{k(j)})\|^2\right) \\
\nonumber
& = {\Lres\gamma \over 2\sqrt{n}{\Lmax}} \sum_{d=k(j)}^{j-1} \E (\|\nabla f(x_{k(d)})\|^2 + \|\nabla f(x_{k(j)})\|^2) \\
\label{eqn_proof_thm_2}
& \leq {\tau \rho^{\tau}\Lres\gamma \over \sqrt{n}{\Lmax}} \E (\|\nabla f(x_{k(j)})\|^2)
\end{align}
where the second line uses \eqref{eqn_proof_thm_1_2}, and the final
inequality uses the fact for $d$ between $k(j)$ and $j-1$, $k(d)$ lies
in the range $k(j)-\tau$ and $j-1$, so we have $|k(d)-k(j)| \le \tau$
for all $d$.

Taking expectation on both sides of \eqref{eqn_proof_thm_1} in terms
of all random variables, together with \eqref{eqn_proof_thm_2}, we
obtain
\begin{align}
\nonumber
&\E (f(x_{j+1}) -f^*) \\
\nonumber
&\quad \leq \E (f(x_j)- f^*) + {\gamma\over n{\Lmax}} \E (|T_1|) - \left({\gamma\over n{\Lmax}} - {\gamma^2\over 2n{\Lmax}}\right)\E (\|\nabla f(x_{k(j)})\|^2)\\
\nonumber
&\quad \leq \E (f(x_j)- f^*) - \left({\gamma\over n{\Lmax}} - {\tau \rho^{\tau}\Lres\gamma^2 \over {n^{3/2}}{\Lmax^2}} - {\gamma^2\over 2n{\Lmax}}\right)\E (\|\nabla f(x_{k(j)})\|^2)\\
\label{eq:fxj1}
&\quad = \E (f(x_j)- f^*)- \frac{\gamma}{n \Lmax} \left( 1-\frac{\psi}{2} \gamma \right)
\E (\|\nabla f(x_{k(j)})\|^2),
\end{align}
which (because of 
\eqnok{eq:boundgamma.1}) implies that $\E (f(x_{j})-f^*)$ is
monotonically decreasing.


Assume now that the essential strong convexity property \eqref{eq:esc}
holds. From Lemma~\ref{lem_1} and the fact that $\E (f(x_{j})-f^*)$ is
monotonically decreasing, we have
\[
\|\nabla f(x_{k(j)})\|^2 \ge
2l\E(f(x_{k(j)})-f^*) \ge
2l\E(f(x_j)-f^*).
\]
So by substituting in \eqref{eq:fxj1}, we obtain
\[
\E (f(x_{j+1}) -f^*) \le \left( 1- \frac{2l\gamma}{n\Lmax} \left(
1-\frac{\psi}{2} \gamma \right) \right) \E (f(x_j)- f^*),
\]
from which the linear convergence claim \eqnok{eqn_thm_2} follows by
an obvious induction.

For the case of general smooth convex $f$, we have from Lemma~\ref{lem_1}, Assumption~\ref{ass_1}, and the monotonic decreasing property for $\E (f(x_{j})-f^*)$ that
\begin{align*}
\E(\| \nabla f(x_{k(j)}) \|^2) \geq & \E\left[\frac{(f(x_{k(j)})-f^*)^2}{\|x_{k(j)}-\PS(x_{k(j)})\|^2}\right] 
\\ \ge &
\frac{\E[(f(x_{k(j)})-f^*)^2]}{R^2} 
\\ \ge &
\frac{[\E(f(x_{k(j)})-f^*)]^2}{R^2} 
\\ \ge &
\frac{[\E(f(x_j)-f^*)]^2}{R^2},
\end{align*}
where the third inequality uses Jensen's inequality $\E(v^2) \geq (\E(v))^2$.
By substituting into \eqref{eq:fxj1}, we obtain
\begin{align*}
\E (f(x_{j+1}) -f^*) 
& \le \E (f(x_j)- f^*)  -
\frac{\gamma}{n\Lmax R^2} \left( 1-\frac{\psi}{2} \gamma \right)
[\E  (f(x_j)-f^*)]^2.
\end{align*}
Defining
\[
C:=\frac{\gamma}{n \Lmax R^2} \left( 1-\frac{\psi}{2} \gamma \right),
\]
we have
\begin{align*}
&\E (f(x_{j+1}) -f^*) \leq \E (f(x_{j}) -f^*) - {C}(\E (f(x_{j}) -f^*))^2\\
\Rightarrow~&{1\over \E (f(x_{j}) -f^*)} \leq {1\over \E (f(x_{j+1}) -f^*)} - {C}{\E (f(x_{j}) -f^*)\over \E (f(x_{j+1}) -f^*)} \\
\Rightarrow~&{1\over \E (f(x_{j+1}) -f^*)} - {1\over \E (f(x_{j}) -f^*)} \geq {C}{\E (f(x_{j}) -f^*)\over \E (f(x_{j+1}) -f^*)} \geq C\\
\Rightarrow~&{1\over \E (f(x_{j+1}) -f^*)} \geq {1\over  f(x_{0}) -f^*} + {C(j+1)} \\
\Rightarrow~&\E (f(x_{j+1}) -f^*) \leq {1\over {(f(x_{0}) -f^*)^{-1}} + {C(j+1)}},
\end{align*}
which completes the proof of the sublinear rate \eqnok{eqn_thm_3}.
\end{proof}

\begin{proof} (Corollary~\ref{co:thm_1})
Note first that for $\rho$ defined by \eqnok{eq:choicerho}, we have
\[
\rho^{\tau} \leq \rho^{\tau+1}=
\left(\left(1+{2e\Lres\over
  \sqrt{n}\Lmax}\right)^{\sqrt{n}\Lmax\over
  2e\Lres}\right)^{{2e\Lres(\tau+1)\over \sqrt{n}\Lmax}} \leq
e^{{2e\Lres(\tau+1)\over \sqrt{n}\Lmax}} \leq e,
\]
and thus from the definition of $\psi$ \eqnok{eq:defpsi} that
\beq \label{eq:psi2}
\psi = 1+ \frac{2 \tau \rho^{\tau} \Lres}{\sqrt{n} \Lmax} \le
1+ \frac{2 \tau e \Lres}{\sqrt{n} \Lmax} \le 2.
\eeq

We show now that the steplength parameter choice $\gamma = 1/\psi$
satisfies all the bounds in \eqnok{eq:boundgamma}, by showing that the
second and third bounds are implied by the first. For 
the second bound \eqnok{eq:boundgamma.2}, we have
\[
\frac{(\rho-1) \sqrt{n} \Lmax}{2 \rho^{\tau+1} \Lres} \ge
\frac{(\rho-1) \sqrt{n} \Lmax}{2 e \Lres} \ge 1 \ge \frac{1}{\psi},
\]
where the second inequality follows from \eqnok{eq:choicerho}.  To
verify that the right hand side of the third bound
\eqnok{eq:boundgamma.3} is greater than $1/\psi$, we consider the
cases $\tau = 0$ and $\tau \geq 1$ separately. For $\tau = 0$, we have
$\psi = 1$ from \eqref{eq:defpsi} and
\[
\frac{(\rho-1) \sqrt{n} \Lmax}{\Lres \rho^{\tau}
  (2+\frac{\Lres}{\sqrt{n} \Lmax})} =
\frac{2e}{2+\frac{\Lres}{\sqrt{n} \Lmax}} \geq \frac{2e}{2+(2e)^{-1}}
\geq 1 \geq \frac{1}{\psi},
\]
where the first inequality is from \eqref{eq:boundtau}.  For the other
case $\tau\geq 1$, we have
\[
\frac{(\rho-1) \sqrt{n} \Lmax}{\Lres \rho^{\tau}
  (2+\frac{\Lres}{\sqrt{n} \Lmax})} = \frac{2e\Lres}{\Lres \rho^{\tau}
  (2+\frac{\Lres}{\sqrt{n} \Lmax})} \ge \frac{2e\Lres}{\Lres e
  (2+\frac{\Lres}{\sqrt{n} \Lmax})} = \frac{2}{2+\frac{\Lres}{\sqrt{n}
    \Lmax}} \ge \frac{1}{\psi}.
\]
We can thus set $\gamma=1/\psi$, and by substituting this choice into
\eqnok{eqn_thm_2} and using \eqnok{eq:psi2}, we obtain
\eqnok{eqn_thm_2_good}.  We obtain \eqnok{eqn_thm_3_good} by making
the same substitution into \eqnok{eqn_thm_3}.
\end{proof}

\begin{proof} (Theorem~\ref{thm:co1})
From  Markov's inequality, we have
\begin{align*}
\mathbb{P}(f(x_j)-f^*\ge \epsilon) & \le \epsilon^{-1}{\E(f(x_j)-f^*)} \\
& \le
{\epsilon}^{-1}{\left(1-{l\over 2n\Lmax}\right)^j(f(x_0)-f^*)}\\
&\le
{\epsilon}^{-1}{\left(1-c\right)^{(1/c) \left|{\log{f(x_0)-f^*\over \eta\epsilon}}\right|}(f(x_0)-f^*)} \quad \mbox{with $c=l/(2n\Lmax)$} \\
& \le
{\epsilon}^{-1}(f(x_0)-f^*) e^{-\left|{\log{f(x_0)-f^*\over \eta\epsilon}}\right|} \\
& =
{\eta}e^{\log{(f(x_0)-f^*)\over \eta\epsilon}} e^{-\left|{\log{f(x_0)-f^*\over \eta\epsilon}}\right|}\\
& \le
\eta,
\end{align*}
where the second inequality applies~\eqref{eqn_thm_2_good}, the third
inequality uses the definition of $j$~\eqref{eq:def_j1}, and the
second last inequality uses the inequality $(1-c)^{1/c} \leq
e^{-1}~\forall \, c\in (0,1)$, which proves the essentially strongly
convex case. Similarly, the general convex case is proven by
\begin{align*}
&\mathbb{P}(f(x_j)-f^*\ge \epsilon) \le \epsilon^{-1}{\E(f(x_j)-f^*)}
\le
\frac{f(x_0)-f^*}{\epsilon\left(1+j{f(x_0)-f^*\over 4n\Lmax R^2}\right)}
\le
\eta,
\end{align*}
where the second inequality uses~\eqref{eqn_thm_3_good} and the last inequality uses the definition of $j$~\eqref{eq:def_j2}.
\end{proof}

\section{Proofs for Constrained Case} \label{app:con}


We start by introducing notation and proving several preliminary
results.  Define
\beq\label{eq:defdelta} 
(\Delta_j)_{i(j)}:=(x_j-x_{j+1})_{i(j)},
\eeq
and formulate the update in Step~\ref{step_proj} of
Algorithm~\ref{alg_ascd} in the following way:
\[
x_{j+1}=\arg\min_{x\in\Omega} \,  \langle \nabla_{i(j)}f(x_{k(j)}), (x-x_j)_{i(j)} \rangle + \frac{\Lmax}{2\gamma}\|x-x_j\|^2.
\]
(Note that $(x_{j+1})_i=(x_j)_i$ for $i\neq i(j)$.) The 
optimality condition for
this formulation is 
\begin{align*}
\left\langle (x-x_{j+1})_{i(j)}, \nabla_{i(j)}f(x_{k(j)}) - {\Lmax\over \gamma}(\Delta_j)_{i(j)} \right\rangle \geq 0, \quad \mbox{\rm for all $x\in \Omega$}.
\end{align*}
This implies in particular that for all $x\in \Omega$, we have
\begin{align}
\left\langle (\PS(x)-x_{j+1})_{i(j)}, \nabla_{i(j)} f(x_{k(j)}) \right\rangle \geq \frac{\Lmax}{\gamma} \left\langle (\PS(x)-x_{j+1})_{i(j)}, (\Delta_j)_{i(j)} \right\rangle.
\label{eqn_pre_1}
\end{align}
From the definition of $\Lmax$, and using the notation
\eqnok{eq:defdelta}, we have
\begin{align*}
f(x_{j+1}) \leq f(x_j) + \langle \nabla_{i(j)}f(x_j), -(\Delta_j)_{i(j)} \rangle + {\Lmax\over 2}\|(\Delta_j)_{i(j)}\|^2,
\end{align*}
which indicates that
\begin{align}
\langle \nabla_{i(j)}f(x_j), (\Delta_j)_{i(j)} \rangle \leq f(x_j) -  f(x_{j+1})  + {\Lmax\over 2}\|(\Delta_j)_{i(j)}\|^2.
\label{eqn_pre_2}
\end{align}
From optimality conditions for the problem \eqref{eq:def.xbar}, which defines
the vector $\bar{x}_{j+1}$, we have 
\begin{align}
\left\langle x-\bar{x}_{j+1},  \nabla f(x_{k(j)}) + \frac{\Lmax}{\gamma}(\bar{x}_{j+1}-x_j) \right\rangle \geq 0\quad \forall \, x\in \Omega.
\label{eqn_pre_3}
\end{align}
We now define $\Delta_j:=x_j-\bar{x}_{j+1}$, and note that this
definition is consistent with $(\Delta)_{i(j)}$ defined
in~\eqref{eq:defdelta}. It can be seen that
\[
\E_{i(j)} (\|x_{j+1}-x_j\|^2) = {1\over n} \|\bar{x}_{j+1}- x_j\|^2.
\]

We now proceed to prove the main results of
Section~\ref{sec_constrained}.

\begin{proof} (Theorem~\ref{thm_2})
We prove \eqref{eqn_thm2_1} by induction. First, note that for any
vectors $a$ and $b$, we have
\begin{align*}
\|a\|^2 - \|b\|^2 = 2\|a\|^2 -(\|a\|^2+\|b\|^2) \leq 2\|a\|^2 -
2\langle a, b\rangle \leq 2\langle a, a-b \rangle \leq 2\|a\|\|a-b\|,
\end{align*}
Thus for all $j$, we have
\begin{equation}
\|x_{j-1}-\bar{x}_j\|^2 - \|x_j-\bar{x}_{j+1}\|^2 \leq 2\|x_{j-1}-\bar{x}_j\| {\|x_j-\bar{x}_{j+1}-x_{j-1}+\bar{x}_j\|}.
\label{eqn_proof2_1}
\end{equation}
The second factor in the r.h.s. of \eqref{eqn_proof2_1} is bounded as follows:
\begin{align}
\nonumber
&\|x_j-\bar{x}_{j+1}-x_{j-1}+\bar{x}_j\|\\
\nonumber
&\quad = \left\|x_j-\PO(x_j-{\gamma \over \Lmax} \nabla f(x_{k(j)}))-(x_{j-1}-\PO(x_{j-1}-{\gamma \over \Lmax} \nabla f(x_{k(j-1)}))) \right\|\\
\nonumber
&\quad \leq \left\|x_j-{\gamma \over \Lmax} \nabla f(x_{k(j)})-\PO(x_j-{\gamma \over \Lmax} \nabla f(x_{k(j)})) -(x_{j-1}-{\gamma \over \Lmax} \nabla f(x_{k(j-1)})\right. \\
\nonumber
&\quad \quad \left.  - \PO(x_{j-1}-{\gamma \over \Lmax} \nabla f(x_{k(j-1)}))) \right\|+ {\gamma \over \Lmax} \left\|\nabla f(x_{k(j-1)}) - \nabla f(x_{k(j)}) \right\| \\
\nonumber
&\quad \leq \left\|x_j-{\gamma \over \Lmax} \nabla f(x_{k(j)})-x_{j-1}+{\gamma \over \Lmax} \nabla f(x_{k(j-1)}) \right\| \\
\nonumber
&\quad\quad+ {\gamma \over \Lmax} \left\| \nabla f(x_{k(j-1)}) - \nabla f(x_{k(j)}) \right\| \\
\nonumber
&\quad \leq \|x_j-x_{j-1}\| + 2{\gamma \over \Lmax} \left\| \nabla f(x_{k(j)}) - \nabla f(x_{k(j-1)}) \right\| \\
\nonumber
&\quad \leq \|x_j-x_{j-1}\| + 2{\gamma \over \Lmax} \sum_{d=\min \{k(j-1), k(j)\}}^{\max \{k(j-1), k(j)\}-1} \|\nabla f(x_{d}) - \nabla f(x_{d+1})\| \\
\label{eqn_proof2_2}
&\quad \leq \|x_j-x_{j-1}\| + 2{\gamma \Lres \over \Lmax}  \sum_{d=\min \{k(j-1), k(j)\}}^{\max \{k(j-1), k(j)\}-1} \|x_{d} -x_{d+1}\|,
\end{align}
where the first inequality follows by adding and subtracting a term,
and the second inequality uses the nonexpansive property of projection:
\[
\|(z-\mathcal{P}_{\Omega}(z)) - (y-\mathcal{P}_{\Omega}(y))\| \leq
\|z-y\|.
\]
One can see that $j-1-\tau\leq k(j-1)\leq j-1$ and $j-\tau\leq
k(j) \leq j$, which implies that $j-1-\tau\leq d\leq j-1$ for each
index $d$ in the summation in \eqnok{eqn_proof2_2}. It also follows that
\begin{align}
 \max \{k(j-1), k(j)\} -1 - \min \{k(j-1), k(j)\} \leq \tau.
\label{eqn_proof2_2_5}
\end{align}

We set $j=1$, and note that $k(0)=0$ and $k(1)\le 1$. Thus, in this
case, we have that the lower and upper limits of the summation in
\eqnok{eqn_proof2_2} are $0$ and $0$, respectively. Thus, this
summation is vacuous, and we have
\[
\|x_1-\bar{x}_{2}+x_{0}-\bar{x}_1\| \leq \left(1+2\frac{\gamma \Lres}{\Lmax}\right)\|x_1-x_{0}\|,
\]
By substituting this bound in \eqref{eqn_proof2_1} and setting $j=1$,
we obtain
\beq \label{eq:crap13}
 \E(\|x_{0}-\bar{x}_1\|^2) - \E(\|x_1-\bar{x}_{2}\|^2) \leq \left(2+4\frac{\gamma \Lres}{\Lmax}\right)\E(\|x_1-x_{0}\|\|\bar{x}_1-x_{0}\|).
\eeq
For any $j$, we have
\begin{align}
\nonumber
\E (\|x_{j}-x_{j-1}\| \|\bar{x}_j -x_{j-1}\|) &\leq {1\over 2}\E(n^{1/2}\|x_{j}-x_{j-1}\|^2 + n^{-1/2}\|\bar{x}_j -x_{j-1}\|^2)\\
\nonumber
& = {1\over 2}\E(n^{1/2}\E_{i(j-1)}(\|x_{j}-x_{j-1}\|^2) + n^{-1/2}\|\bar{x}_j -x_{j-1}\|^2)\\
\nonumber
& = {1\over 2}\E(n^{-1/2}\|\bar{x}_{j}-x_{j-1}\|^2 + n^{-1/2}\|\bar{x}_j -x_{j-1}\|^2)\\
\label{eqn_proof2_4}
& = n^{-1/2}\E \|\bar{x}_j -x_{j-1}\|^2.
\end{align}
Returning to \eqnok{eq:crap13}, we have
\[
 \E(\|x_{0}-\bar{x}_1\|^2) - \E(\|x_1-\bar{x}_{2}\|^2) \leq \left({2\over \sqrt{n}}+\frac{4\gamma\Lres}{\sqrt{n}\Lmax}\right)\E \|\bar{x}_1 -x_{0}\|^2
\]
which implies that
\[
 \E(\|x_{0}-\bar{x}_1\|^2) \leq \left(1- {2\over \sqrt{n}}-\frac{4\gamma\Lres}{\sqrt{n}\Lmax}\right)^{-1} \E(\|x_1-\bar{x}_{2}\|^2) \leq \rho  \E(\|x_1-\bar{x}_{2}\|^2).
\]
To see the last inequality above, we only need to verify that
\[
\gamma \leq \left(1-\rho^{-1}-{2\over \sqrt{n}}\right)\frac{\sqrt{n}\Lmax}{4\Lres}.
\]
This proves that \eqref{eqn_thm2_1} holds for $j=1$.

To take the inductive step, we assume that
\eqref{eqn_thm2_1} holds up to index $j-1$. We have for $j-1-\tau \leq
d \leq j-2$ that
\begin{align}
\nonumber
\E (\|x_{d}-x_{d+1}\| \|\bar{x}_j -x_{j-1}\|)
& \leq {1\over 2}\E (n^{1/2}\|x_{d}-x_{d+1}\|^2 + n^{-1/2}\|\bar{x}_j -x_{j-1}\|^2)\\
\nonumber
& = {1\over 2} \E (n^{1/2}\E_{i(d)}(\|x_{d}-{x}_{d+1}\|^2) + n^{-1/2}\|\bar{x}_j -x_{j-1}\|^2) \\
\nonumber
& = {1\over 2} \E (n^{-1/2}\|x_{d}-\bar{x}_{d+1}\|^2 + n^{-1/2}\|\bar{x}_j -x_{j-1}\|^2) \\
\nonumber
& \leq {1\over 2} \E (n^{-1/2}\rho^{\tau}\|x_{j-1}-\bar{x}_{j}\|^2 + n^{-1/2}\|\bar{x}_j -x_{j-1}\|^2)\\
& \leq {\rho^{\tau}\over n^{1/2}} \E (\|\bar{x}_j -x_{j-1}\|^2),
\label{eqn_proof2_3}
\end{align}
where the second inequality uses the inductive hypothesis.
By substituting \eqref{eqn_proof2_2} into \eqref{eqn_proof2_1} and
taking expectation on both sides of~\eqref{eqn_proof2_1}, we obtain
\begin{align*}
&\E(\|x_{j-1}-\bar{x}_j\|^2) - \E(\|x_j-\bar{x}_{j+1}\|^2)\\
&\quad \leq 2\E(\|\bar{x}_j -x_{j-1}\| \|\bar{x}_j-\bar{x}_{j+1}+x_j-x_{j-1}\|)\\
&\quad \leq 2\E\left(\|\bar{x}_j -x_{j-1}\| \left(\|x_j-x_{j-1}\| + 2{\gamma \Lres \over \Lmax} \sum_{d=\min \{k(j-1), k(j)\}}^{\max \{k(j-1), k(j)\}-1} \|x_{d} -x_{d+1}\|\right)\right)\\
&\quad = 2\E(\|\bar{x}_j -x_{j-1}\| \|x_j-x_{j-1}\|) + \\
\nonumber
&\quad\quad\quad\quad 4{\gamma \Lres \over \Lmax} \E\left(\sum_{d=\min \{k(j-1), k(j)\}}^{\max \{k(j-1), k(j)\}-1} (\|\bar{x}_j -x_{j-1}\|\|x_{d} -x_{d+1}\|)\right) \\
&\quad \leq n^{-1/2}\left(2+\frac{4\gamma \Lres\tau \rho^{\tau}}{\Lmax}\right) \E (\|{x}_{j-1}-\bar{x}_j\|^2),
\end{align*}
where the last line uses \eqref{eqn_proof2_2_5}, \eqref{eqn_proof2_4}, and \eqref{eqn_proof2_3}. It follows that
\begin{align*}
&\E(\|x_{j-1}-\bar{x}_j\|^2) \leq \left(1-n^{-1/2}\left(2+\frac{4\gamma
  \Lres\tau \rho^{\tau}}{\Lmax}\right)\right)^{-1} \E
(\|x_j-\bar{x}_{j+1}\|^2) \leq \rho \E (\|x_j-\bar{x}_{j+1}\|^2).
\end{align*}
To see the last inequality, one only needs to verify that
\[
\rho^{-1} \leq 1-{1\over \sqrt{n}}\left(2+{4\gamma \Lres\tau\rho^{\tau}\over \Lmax}\right) \;
\Leftrightarrow \;  \gamma \leq \left({1-\rho^{-1} - {2\over \sqrt{n}}}\right){{\sqrt{n}\Lmax \over 4\Lres\tau\rho^{\tau}}},
\]
and the last inequality is true because of the upper bound of $\gamma$ in \eqref{eq:boundgammac}. It proves \eqref{eqn_thm2_1}.

Next we will show the expectation of objective is monotonically
decreasing. We have using the definition \eqnok{eq:defdelta} that
\begin{align}
\nonumber
&\E_{i(j)}(f(x_{j+1})) = n^{-1}\sum_{i=1}^n f(x_j + (\Delta_j)_i)\\
\nonumber
&\quad \leq n^{-1}\sum_{i=1}^n\left[f(x_j) + \langle \nabla_{i}f(x_j), (\bar{x}_{j+1}-x_j)_{i} \rangle + \frac{\Lmax}{2}\|(x_{j+1}-x_j)_{i}\|^2\right] \\
\nonumber
&\quad = f(x_j)+n^{-1} \left( \langle \nabla f(x_j), \bar{x}_{j+1}-x_j \rangle + \frac{\Lmax}{2}\|\bar{x}_{j+1}-x_j\|^2 \right) \\
\nonumber
&\quad = f(x_j) + {1\over n}\left(\langle \nabla f(x_{k(j)}), \bar{x}_{j+1}-x_j \rangle + \frac{\Lmax}{2}\|\bar{x}_{j+1}-x_j\|^2 \right) +{1\over n} \langle \nabla f(x_j)- \nabla f(x_{k(j)}), \bar{x}_{j+1}-x_j\rangle \\
\nonumber
&\quad \leq f(x_j) + {1\over n}\left(\frac{\Lmax}{2}\|\bar{x}_{j+1}-x_j\|^2  - {\Lmax\over \gamma}\|\bar{x}_{j+1}-x_j\|^2\right) 
+ {1\over n} \langle \nabla f(x_j)- \nabla f(x_{k(j)}), \bar{x}_{j+1}-x_j\rangle\\
&\quad = f(x_j) -  \left({\frac{1}{\gamma}-\frac{1}{2}}\right)\frac{\Lmax}{n}\|\bar{x}_{j+1}-x_j\|^2 +  {1\over n}\langle \nabla f(x_j)- \nabla f(x_{k(j)}), \bar{x}_{j+1}-x_j\rangle,
\label{eqn_proof2_5}
\end{align}
where the second inequality uses~\eqref{eqn_pre_3}. Consider the expectation of the last term on the right-hand side of this expression. We have
\begin{align}
\nonumber
&\E  \langle \nabla f(x_j)- \nabla f(x_{k(j)}), \bar{x}_{j+1}-x_j\rangle\\
\nonumber
&\quad \leq \E \left(\| \nabla f(x_j)- \nabla f(x_{k(j)})\| \|\bar{x}_{j+1}-x_j \|\right) \\
\nonumber
&\quad \leq \E \left(\sum_{d=k(j)}^{j-1}\|\nabla f(x_d) - \nabla f(x_{d+1})\| \|\bar{x}_{j+1}-x_j\|\right)\\
\nonumber
&\quad \leq \Lres\E \left(\sum_{d=k(j)}^{j-1} \|x_d-x_{d+1}\| \|\bar{x}_{j+1}-x_j\|\right)\\
\nonumber
&\quad \leq {\Lres\over 2} \E \left(\sum_{d=k(j)}^{j-1} (n^{1/2}\|x_d-x_{d+1}\|^2 + n^{-1/2}\|\bar{x}_{j+1}-x_j\|^2)\right) \\
\nonumber
&\quad = {\Lres\over 2} \E \left(\sum_{d=k(j)}^{j-1} \left(n^{1/2}\E_{i(d)}(\|x_d-x_{d+1}\|^2 ) + n^{-1/2}\|\bar{x}_{j+1}-x_j\|^2 \right)\right) \\
\nonumber
&\quad = {\Lres\over 2} \E \left(\sum_{d=k(j)}^{j-1} (n^{-1/2}\|x_d-\bar{x}_{d+1}\|^2 + n^{-1/2}\|\bar{x}_{j+1}-x_j\|^2)\right) \\
\nonumber
&\quad \leq {\Lres\over 2n^{1/2}} \E \left(\sum_{d=k(j)}^{j-1} (1+\rho^{\tau}) \|\bar{x}_{j+1}-x_j\|^2 \right)\\
&\quad \leq {\Lres\tau\rho^{\tau}\over n^{1/2}} \E \|\bar{x}_{j+1}-x_j\|^2,
\label{eqn_proof2_5_5}
\end{align}
where the fifth inequality uses \eqref{eqn_thm2_1}. By taking
expectation on both sides of \eqref{eqn_proof2_5} and substituting
\eqref{eqn_proof2_5_5}, we have
\begin{align*}
&\E(f(x_{j+1})) \leq \E (f(x_{j})) - {1\over n}\left(\left({\frac{1}{\gamma}-\frac{1}{2}}\right)\Lmax-{\Lres\tau\rho^{\tau}\over n^{1/2}}\right)\E\|\bar{x}_{j+1}-x_j\|^2.
\end{align*}
To see $\left({\frac{1}{\gamma}-\frac{1}{2}}\right)\Lmax-{\Lres\tau\rho^{\tau}\over n^{1/2}}\geq 0$, we only need to verify
\begin{align*}
&\gamma \leq \left({1\over 2} + {\Lres\tau\rho^{\tau}\over \sqrt{n}\Lmax}\right)^{-1} \nonumber
\end{align*}
which is implied by the first upper bound of
$\gamma$~\eqref{eq:boundgammac}.  Therefore, we have proved the
monotonicity $\E(f(x_{j+1})) \leq \E (f(x_{j}))$.

Next we prove the sublinear convergence rate for the constrained
smooth convex case in~\eqref{eqn_thm2_2}. We have
\begin{align}
\nonumber
& \|x_{j+1} - \PS(x_{j+1})\|^2 \\
\nonumber &\quad 
 \leq \|x_{j+1} - \PS(x_j)\|^2 \\
\nonumber
&\quad = \|x_j-(\Delta_j)_{i(j)}e_{i(j)} - \PS(x_j)\|^2 \\
\nonumber
&\quad = \|x_j - \PS(x_j)\|^2 + |(\Delta_j)_{i(j)}|^2 - 2  (x_j - \PS(x_j))_{i(j)} (\Delta_j)_{i(j)} \\
\nonumber
&\quad = \|x_j - \PS(x_j)\|^2 - |(\Delta_j)_{i(j)}|^2 - 2  \left( (x_j - \PS(x_j))_{i(j)}-(\Delta_j)_{i(j)} \right) (\Delta_j)_{i(j)} \\
\nonumber
&\quad = \|x_j - \PS(x_j)\|^2 - |(\Delta_j)_{i(j)}|^2 + 2  (\PS(x_j)-x_{j+1})_{i(j)} (\Delta_j)_{i(j)} \\
&\quad \leq \|x_j - \PS(x_j)\|^2 - |(\Delta_j)_{i(j)}|^2 + \frac{2\gamma}{\Lmax} (\PS(x_j)-x_{j+1})_{i(j)} \nabla_{i(j)} f(x_{k(j)})
\label{eq_proof2_6a}
\end{align}
where the last inequality uses \eqref{eqn_pre_1}.
We bound the last term in \eqref{eq_proof2_6a} by
\begin{align}
\nonumber
& \frac{2\gamma}{\Lmax} (\PS(x_j)-x_{j+1})_{i(j)} \nabla_{i(j)} f(x_{k(j)}) 
\\ \nonumber & \quad =
\frac{2\gamma}{\Lmax} (\PS(x_j)-x_j)_{i(j)} \nabla_{i(j)} f(x_{k(j)})  \\
\nonumber
&\quad\quad + \frac{2\gamma}{\Lmax}\left((\Delta_j)_{i(j)} \nabla_{i(j)} f(x_{j})  +
(\Delta_j)_{i(j)} \left(\nabla_{i(j)} f(x_{k(j)}) - \nabla_{i(j)} f(x_{j}) \right) \right) \\
\nonumber & \quad \leq 
 \frac{2\gamma}{\Lmax}  (\PS(x_j)-x_j)_{i(j)} \nabla_{i(j)} f(x_{k(j)}) \\
\nonumber
&\quad\quad + \frac{2\gamma}{\Lmax}\bigg(f(x_j) -  f(x_{j+1})  + {\Lmax\over 2}|(\Delta_j)_{i(j)}|^2 \\
\nonumber
&\quad\quad +  (\Delta_j)_{i(j)} \left( \nabla_{i(j)} f(x_{k(j)}) - \nabla_{i(j)} f(x_{j}) \right) \bigg)
\\ \nonumber
&\quad =
\gamma|(\Delta_j)_{i(j)}|^2  +\frac{2\gamma}{\Lmax}(f(x_j) -  f(x_{j+1}) ) 
\\ \nonumber &\quad \quad
+ \frac{2\gamma}{\Lmax} (\PS(x_j)-x_j)_{i(j)} \nabla_{i(j)} f(x_{k(j)})
\\ &\quad\quad
+ \frac{2\gamma}{\Lmax}(\Delta_j)_{i(j)} \left( \nabla_{i(j)} f(x_{k(j)}) - \nabla_{i(j)} f(x_{j}) \right) 
\label{eq_proof2_6b}
\end{align}
where the inequality uses \eqref{eqn_pre_2}. 

Together with \eqref{eq_proof2_6a}, we obtain
\begin{align}
\nonumber
&\|x_{j+1} - \PS(x_{j+1})\|^2 \leq \|x_j - \PS(x_j)\|^2 - (1-\gamma)|(\Delta_j)_{i(j)}|^2 
\\ \nonumber
&\quad\quad  +\frac{2\gamma}{\Lmax}(f(x_j) -  f(x_{j+1}) ) 
\\ &\quad \quad
+ \frac{2\gamma}{\Lmax} \underbrace{(\PS(x_j)-x_j)_{i(j)} \nabla_{i(j)} f(x_{k(j)})}_{T_1}
\label{eqn_proof2_6}
\\ \nonumber &\quad\quad
+ \frac{2\gamma}{\Lmax}\underbrace{ (\Delta_j)_{i(j)} \left( \nabla_{i(j)} f(x_{k(j)}) - \nabla_{i(j)} f(x_{j}) \right)}_{T_2}.
\end{align}
We now seek upper bounds on the
quantities $T_1$ and $T_2$ in the expectation sense. For $T_1$, we have

\begin{align*}
& \E(T_1) 
= n^{-1}\E \langle \PS(x_j)-x_j, \nabla f(x_{k(j)})  \rangle \\
& = n^{-1}\E \langle \PS(x_j)-x_{k(j)}, \nabla f(x_{k(j)})  \rangle  + n^{-1}\E \sum_{d={k(j)}}^{j-1}\langle x_{d}-x_{d+1}, \nabla f(x_{k(j)})  \rangle \\
& = n^{-1}\E \langle \PS(x_j)-x_{k(j)}, \nabla f(x_{k(j)})  \rangle  \\
& \quad + n^{-1}\E \sum_{d={k(j)}}^{j-1} \langle x_{d}-x_{d+1}, \nabla f(x_{d})  \rangle+ \langle x_{d}-x_{d+1}, \nabla f(x_{k(j)}) - \nabla f(x_d)  \rangle \\
& \leq n^{-1}\E (f^*- f(x_{k(j)})) + n^{-1} \E \sum_{d={k(j)}}^{j-1} \left( f(x_{d})-f(x_{d+1}) + {\Lmax\over 2}\|x_d-x_{d+1}\|^2 \right) \\
& \quad +n^{-1} \E \sum_{d={k(j)}}^{j-1}  \langle x_{d}-x_{d+1}, \nabla f(x_{k(j)}) - \nabla f(x_d)  \rangle\\
&= n^{-1}\E (f^*- f(x_{j})) + {\Lmax\over 2n} \E \sum_{d={k(j)}}^{j-1}\|x_d-x_{d+1}\|^2 
\\ &\quad 
+n^{-1} \E \sum_{d={k(j)}}^{j-1}  \langle x_{d}-x_{d+1}, \nabla f(x_{k(j)}) - \nabla f(x_d)  \rangle\\
& = n^{-1}\E (f^*- f(x_{j})) + {\Lmax\over 2n^2} \E \sum_{d={k(j)}}^{j-1}\|x_d-\bar{x}_{d+1}\|^2 
\\ &\quad 
+n^{-1} \E \sum_{d={k(j)}}^{j-1}  \langle x_{d}-x_{d+1}, \nabla f(x_{k(j)}) - \nabla f(x_d)  \rangle\\
& \leq n^{-1}\E (f^*- f(x_{j})) + {\Lmax\tau\rho^{\tau}\over 2n^2} \E \|x_j-\bar{x}_{j+1}\|^2 
\\ &\quad 
+n^{-1} \sum_{d={k(j)}}^{j-1}  \underbrace{\E \langle x_{d}-x_{d+1}, \nabla f(x_{k(j)}) - \nabla f(x_d)  \rangle}_{T_3},
\end{align*}
where the first inequality uses the convexity of $f(x)$: 
\[
f(\PS(x_j)) - f(x_{k_{j}}) \geq  \langle \PS(x_j)-x_{k(j)}, \nabla f(x_{k(j)})  \rangle,
\]
and the last inequality uses \eqref{eqn_thm2_1}. 

The upper bound of $\E(T_3)$ is estimated by
\begin{align*}
\E(T_3) & = \E  \langle x_{d}-x_{d+1}, \nabla f(x_{k(j)}) - \nabla f(x_d)  \rangle \\
& = \E(\E_{i(d)} \langle x_{d}-x_{d+1}, \nabla f(x_{k(j)}) - \nabla f(x_d)  \rangle) \\
& = n^{-1}\E \langle x_{d}-\bar{x}_{d+1}, \nabla f(x_{k(j)}) - \nabla f(x_d)  \rangle \\
& \leq n^{-1}\E \|x_{d}-\bar{x}_{d+1}\| \|\nabla f(x_{k(j)}) - \nabla f(x_d)\| \\
& \leq n^{-1} \E\left(\|x_{d}-\bar{x}_{d+1}\| \sum_{t=k(j)}^{d-1}\|\nabla f(x_{t}) - \nabla f(x_{t+1})\|\right)\\
& \leq \frac{\Lres}{n} \E\left(\sum_{t=k(j)}^{d-1}\|x_{d}-\bar{x}_{d+1}\| \|x_{t} - x_{t+1}\|\right) \\
& \leq \frac{\Lres}{2n} \E\left(\sum_{t=k(j)}^{d-1}(n^{-1/2}\|x_{d}-\bar{x}_{d+1}\|^2 + n^{1/2} \|x_{t} - x_{t+1}\|^2)\right) \\
& \leq \frac{\Lres}{2n} \sum_{t=k(j)}^{d-1}\E(n^{-1/2}\|x_{d}-\bar{x}_{d+1}\|^2 + n^{-1/2} \|x_{t} - \bar{x}_{t+1}\|^2) \\
& \leq \frac{\Lres\rho^{\tau}}{n^{3/2}} \sum_{t=k(j)}^{d-1}\E(\|x_{j}-\bar{x}_{j+1}\|^2) \\
& \leq \frac{\Lres\tau\rho^{\tau}}{n^{3/2}} \E(\|x_{j}-\bar{x}_{j+1}\|^2),
\end{align*}
where the first inequality uses the Cauchy inequality and the second last inequality uses \eqref{eqn_thm2_1}. Therefore, $\E(T_1)$ can be bounded by
\begin{align}
\nonumber
\E(T_1) & = \E\langle (\PS(x_j)-x_j)_{i(j)}, \nabla_{i(j)} f(x_{k(j)}) \rangle \\
\nonumber
& \leq \frac{1}{n}\E (f^*- f(x_{j})) + {\Lmax\tau\rho^{\tau}\over 2n^2} \E \|x_j-\bar{x}_{j+1}\|^2 
+ \sum_{d=k(j)}^{j-1}\frac{\Lres\tau\rho^{\tau}}{n^{5/2}} \E(\|x_{j}-\bar{x}_{j+1}\|^2) \\
\label{eq:bT1}
& \leq \frac{1}{n}\left(f^*-\E f(x_{j}) + \left({\Lmax\tau\rho^{\tau}\over 2n} + \frac{\Lres\tau^2\rho^{\tau}}{n^{3/2}}\right) \E(\|x_{j}-\bar{x}_{j+1}\|^2)\right).
\end{align}

For $T_2$, we have
\begin{align}
\nonumber
\E(T_2) & =\E  (\Delta_j)_{i(j)} \left( \nabla_{i(j)} f(x_{k(j)}) - \nabla_{i(j)} f(x_{j}) \right) \\
\nonumber
& = n^{-1}\E \langle \Delta_j, \nabla f(x_{k(j)}) - \nabla f(x_j) \rangle \\
\nonumber
& \leq n^{-1} \E (\|\Delta_j\| \|\nabla f(x_{k(j)}) - \nabla f(x_j)\|) \\
\nonumber
& \leq n^{-1} \E \left( \sum _{d=k(j)}^{j-1} \|\Delta_j\| \|\nabla f(x_{d}) - \nabla f(x_{d+1})\| \right)\\
\nonumber
& \leq \frac{\Lres}{n} \E \left( \sum _{d=k(j)}^{j-1}  \|\Delta_j\| \|x_{d} - x_{d+1}\| \right) \\
\nonumber
& = \frac{\Lres}{2n} \E \left( \sum _{d=k(j)}^{j-1} n^{-1/2}\|\Delta_j\|^2 + n^{1/2}\|x_{d} - x_{d+1}\|^2 \right) \\
\nonumber
& = \frac{\Lres}{2n} \E \left( \sum _{d=k(j)}^{j-1} n^{-1/2}\|x_j-\bar{x}_{j+1}\|^2 + n^{1/2}\E_{i(d)}\|x_{d} - x_{d+1}\|^2 \right) \\
\nonumber
& = \frac{\Lres}{2n} \E \left( \sum _{d=k(j)}^{j-1} n^{-1/2}\|x_j-\bar{x}_{j+1}\|^2 + n^{-1/2}\|x_{d} - \bar{x}_{d+1}\|^2 \right) \\
\nonumber
& = \frac{\Lres}{2n^{3/2}} \left( \sum _{d=k(j)}^{j-1} \E\|x_j-\bar{x}_{j+1}\|^2 + \E\|x_{d} - \bar{x}_{d+1}\|^2 \right) \\
\nonumber
& \leq \frac{\Lres(1+\rho^{\tau})}{2n^{3/2}} \sum _{d=k(j)}^{j-1} \E\|x_j-\bar{x}_{j+1}\|^2 \\
\label{eq:bT2}
& \leq \frac{\Lres\tau\rho^{\tau}}{n^{3/2}} \E\|x_j-\bar{x}_{j+1}\|^2,
\end{align}
where the second last inequality uses \eqref{eqn_thm2_1}.

%
%

By taking the expectation on both sides of \eqref{eqn_proof2_6}, using
$\E_{i(j)}(|(\Delta_j)_{i(j)}|^2) = n^{-1} \| x_j-\bar{x}_{j+1} \|^2$,
and substituting the upper bounds from \eqnok{eq:bT1} and
\eqnok{eq:bT2}, we obtain
\begin{align}
\nonumber
\E\|x_{j+1} & - \PS(x_{j+1})\|^2 \leq \E \|x_j - \PS(x_j)\|^2 \\
\nonumber
& \quad - {1\over n}\left(1-\gamma - \frac{2\gamma \Lres\tau\rho^{\tau}}{\Lmax n^{1/2}} -{{\gamma \tau\rho^{\tau}\over n} - \frac{2\gamma \Lres\tau^2\rho^{\tau}}{\Lmax n^{3/2}}} \right)\E\|x_j-\bar{x}_{j+1}\|^2 \\
\nonumber
&\quad +\frac{2\gamma}{\Lmax n}(f^*- \E f(x_j)) + \frac{2\gamma}{\Lmax}(\E f(x_j) -  \E f(x_{j+1}) )  \\
\label{eqn_proof2_7_5}
&\leq \E \|x_j - \PS(x_j)\|^2 + \frac{2\gamma}{\Lmax n}(f^*- \E f(x_j)) + \frac{2\gamma}{\Lmax}(\E f(x_j) -  \E f(x_{j+1}) ).
\end{align}
In the second inequality, we were able to drop the term involving $\E
\|x_j-\bar{x}_{j+1}\|^2$ by using the fact that
\[
1-\gamma - \frac{2\gamma \Lres\tau\rho^{\tau}}{\Lmax n^{1/2}} -{{\gamma
    \tau\rho^{\tau}\over n} - \frac{2\gamma
    \Lres\tau^2\rho^{\tau}}{\Lmax n^{3/2}}} = 1 - \gamma\psi \geq 0,
\]
which follows from the definition \eqnok{eq:defpsic} of $\psi$ and
from the first upper bound on $\gamma$ in \eqref{eq:boundgammac}.  It
follows that
\begin{align}
\nonumber
&\E\|x_{j+1} - \PS(x_{j+1})\|^2 + \frac{2\gamma}{\Lmax}(\E f(x_{j+1}) -f^*) \\
\label{eq:crap6}
&\quad \leq  \E \|x_j - \PS(x_j)\|^2 + \frac{2\gamma}{\Lmax}(\E f(x_j) - f^*) - \frac{2\gamma}{\Lmax n}( \E f(x_j)- f^*) \\
\nonumber
&\quad \leq \|x_0 - \PS(x_0)\|^2 +  \frac{2\gamma}{\Lmax}(f(x_0) - f^*) - \frac{2\gamma}{\Lmax n}\sum_{t=0}^{j}( \E f(x_t)- f^*)\\
\nonumber
&\quad \leq R_0^2 + \frac{2\gamma}{\Lmax}(f(x_0) - f^*) - \frac{2\gamma (j+1)}{\Lmax n}(\E f(x_{j+1})- f^*),
\end{align}
where the second inequality follows by applying induction to the
inequality
\[
S_{j+1}\leq S_j-{2\gamma\over \Lmax n}\E (f(x_j)-f^*),
\]
where
\[
S_j:=\E (\|x_j-\PS(x_j)\|^2) + {2\gamma \over \Lmax}\E
(f(x_j)-\PS(x_j)),
\]
and the last line uses the monotonicity of $\E f(x_{j})$ (proved
above) and the definition $R_0 = \|x_0-\PS(x_0)\|$. It implies that
\begin{alignat*}{2}
&\quad\E\|x_{j+1} - \PS(x_{j+1})\|^2 + \frac{2\gamma}{\Lmax}(\E f(x_{j+1}) -f^*) + \frac{2\gamma (j+1)}{\Lmax n}(\E f(x_{j+1})- f^*)\\
& \quad \quad \leq R_0^2 + \frac{2\gamma}{\Lmax}(f(x_0) - f^*) \\
&\Rightarrow  \frac{2\gamma (n+j+1)}{\Lmax n}(\E f(x_{j+1})- f^*) \leq R_0^2 + \frac{2\gamma}{\Lmax}(f(x_0) - f^*) \\
& \Rightarrow \E f(x_{j+1})- f^* \leq \frac{n(R_0^2\Lmax+ 2\gamma(f(x_0)- f^*))}{2\gamma (n+j+1)}.
\end{alignat*}
This completes the proof of  the sublinear convergence rate \eqnok{eqn_thm2_2}.

Finally, we prove the linear convergence rate \eqnok{eqn_thm2_3} for
the essentially strongly convex case. All bounds proven above hold,
and we make use of the following additional property:
\begin{align*}
f(x_j) - f^* \geq \langle \nabla f(\PS(x_j)), x_j - \PS(x_j)\rangle + \frac{l}{2}\|x_j - \PS(x_j)\|^2 \geq \frac{l}{2}\|x_j - \PS(x_j)\|^2,
\end{align*}
due to feasibility of $x_j$ and $\langle \nabla f(\PS(x_j)), x_j -
\PS(x_j)\rangle \geq 0$. By using this result together with some
elementary manipulation, we obtain
\begin{align}
\nonumber
f(x_j) - f^* &= \left(1- \frac{\Lmax}{l\gamma+\Lmax}\right)(f(x_j) - f^*) +   \frac{\Lmax}{l\gamma+\Lmax}(f(x_j) - f^*)\\
\nonumber
&\geq \left(1- \frac{\Lmax}{l\gamma+\Lmax}\right)(f(x_j) - f^*) +   \frac{\Lmax l}{2(l\gamma+\Lmax)}\|x_j - \PS(x_j)\|^2\\
\label{eqn_proof2_9}
& = \frac{\Lmax l}{2(l\gamma+\Lmax)}\left(\|x_j - \PS(x_j)\|^2 + \frac{2\gamma}{\Lmax}(f(x_j) - f^*)\right).
\end{align}
Recalling \eqnok{eq:crap6}, we have
\begin{align}
\nonumber
&\E\|x_{j+1} - \PS(x_{j+1})\|^2 + \frac{2\gamma}{\Lmax}(\E f(x_{j+1}) -f^*) \\
\label{eq:crap7}
&\quad \leq  \E \|x_j - \PS(x_j)\|^2 + \frac{2\gamma}{\Lmax}(\E f(x_j) - f^*) - \frac{2\gamma}{\Lmax n}( \E f(x_j)- f^*).
\end{align}
By taking the expectation of both sides in \eqref{eqn_proof2_9} and
substituting in the last term of \eqnok{eq:crap7}, we obtain
\begin{align*}
&\E\|x_{j+1} - \PS(x_{j+1})\|^2 + \frac{2\gamma}{\Lmax}(\E f(x_{j+1}) -f^*)\\
&\quad \leq  \E \|x_j - \PS(x_j)\|^2 + \frac{2\gamma}{\Lmax}(\E f(x_j) - f^*)\\
&\quad \quad - \frac{2\gamma}{\Lmax n}\left(\frac{\Lmax l}{2(l\gamma+\Lmax)}\left(\E\|x_j - \PS(x_j)\|^2 + \frac{2\gamma}{\Lmax}(\E f(x_j) - f^*)\right)\right)\\
&\quad = \left(1- \frac{l} {n(l+\gamma^{-1}\Lmax)} \right) \left(\E \|x_j - \PS(x_j)\|^2 + \frac{2\gamma}{\Lmax}(\E f(x_j) - f^*)\right) \\
&\quad \leq \left(1- \frac{l} {n(l+\gamma^{-1}\Lmax)} \right)^{j+1} \left(\|x_0-\PS(x_0)\|^2 + \frac{2\gamma}{\Lmax}(f(x_0) - f^*)\right),
\end{align*}
which yields \eqref{eqn_thm2_3}.
\end{proof}

\begin{proof} (Corollary~\ref{co:thm_2})
To apply Theorem~\ref{thm_2}, we first show $\rho > \left(1-{2\over
  \sqrt{n}}\right)^{-1}$. Using the bound \eqnok{eq:boundtauc},
together with $\Lres/\Lmax \ge 1$, we obtain
\begin{align*}
&\left(1-{2\over \sqrt{n}}\right)\left(1+{4e\tau \Lres\over  \sqrt{n}\Lmax}\right)
=
\left(1+{4e\tau \Lres\over  \sqrt{n}\Lmax}\right) - \left(1+{4e\tau \Lres\over  \sqrt{n}\Lmax}\right){2\over \sqrt{n}}
\\ &\quad \ge
\left(1+{4e\tau \over  \sqrt{n}}\right) - \left(1+{1\over  \tau+1}\right){2\over \sqrt{n}} =1+ \left({2e\tau } - 1 -{1\over \tau+1}\right){2\over \sqrt{n}} > 1,
\end{align*}
where the last inequality uses $\tau\ge 1$. Note that for $\rho$
defined by~\eqref{eq:choicerhoc}, and using \eqnok{eq:boundtauc}, we
have
\[
\rho^{\tau} \leq \rho^{\tau+1}= \left(\left(1+{4e\tau \Lres\over
  \sqrt{n}\Lmax}\right)^{\sqrt{n}\Lmax\over 4e\tau
  \Lres}\right)^{{4e\tau \Lres(\tau+1)\over \sqrt{n}\Lmax}} \leq e^{{4e\tau
    \Lres(\tau+1)\over \sqrt{n}\Lmax}} \leq e.
\]
Thus from the definition of $\psi$~\eqref{eq:defpsic},  we have that
\begin{align}
\nonumber
&\psi = {1+\frac{\Lres\tau\rho^{\tau}}{\sqrt{n}\Lmax}\left({2 +{\Lmax\over \sqrt{n}\Lres} + {2\tau\over n}}\right)}
\leq
1+ \frac{\Lres\tau\rho^{\tau}}{4e\Lres\tau(\tau+1)}\left(2+{1\over {\sqrt{n}}} + {2\tau\over n}\right)
\\
\label{eq:psi2c}
&\quad \leq
1+ \frac{1}{4(\tau+1)}\left(2+{1\over {\sqrt{n}}} + {2\tau\over n}\right)
\leq 1 + \left({1\over 4} + {1\over 16} + {1\over 10}\right) \leq 2.
\end{align}
(The second last inequality uses $n\ge 5$ and $\tau\ge 1$.)
Thus, the steplength parameter choice $\gamma=1/2$ satisfies the first
bound in~\eqref{eq:boundgammac}. To show that the second bound in
\eqnok{eq:boundgammac} holds also, we have
\begin{align*}
&\left(1-{1\over \rho}-{2\over \sqrt{n}}\right)\frac{\sqrt{n}\Lmax}{4\Lres\tau\rho^{\tau}} = \left({\rho-1\over \rho}-{2\over \sqrt{n}}\right)\frac{\sqrt{n}\Lmax}{4\Lres\tau\rho^{\tau}}
\\
&\quad =
\frac{4e\tau \Lres}{4\Lres\tau \rho^{\tau+1}} - \frac{\Lmax}{2\Lres\tau\rho^{\tau}}
\ge
1 - {1\over 2} ={1\over 2}.
\end{align*}
We can thus set $\gamma=1/2$, and by substituting this choice
into~\eqref{eqn_thm2_3}, we obtain~\eqref{eqn_thm_2_good_c}. We
obtain~\eqref{eqn_thm_3_good_c} by making the same substitution
into~\eqref{eqn_thm2_2}.
\end{proof}

\vskip 0.2in
{
\bibliographystyle{plainnat}
\bibliography{reference}
}

\end{document}